\documentclass[11pt]{article}
\usepackage{definitions}

\title{cuPDLP.jl: A GPU Implementation of Restarted Primal-Dual Hybrid Gradient for Linear Programming in Julia}
\author{Haihao Lu\thanks{The University of Chicago, Booth School of Business (haihao.lu@chicagobooth.edu).} \and Jinwen Yang\thanks{The University of Chicago, Department of Statistics (jinweny@uchicago.edu).}}
\date{June 2024}

\begin{document}

\maketitle

\begin{abstract}
    In this paper, we provide an affirmative answer to the long-standing question: \emph{Are GPUs useful in solving linear programming?} We present cuPDLP.jl, a GPU implementation of restarted primal-dual hybrid gradient (PDHG) for solving linear programming (LP). We show that this prototype implementation in Julia has comparable numerical performance on standard LP benchmark sets to Gurobi, a highly optimized implementation of the simplex and interior-point methods. This demonstrates the power of using GPUs in linear programming, which, for the first time, showcases that GPUs and first-order methods can lead to performance comparable to state-of-the-art commercial optimization LP solvers on standard benchmark sets.
\end{abstract}

\section{Introduction}
Linear programming (LP) is a fundamental optimization problem class with a long history and a vast range of applications in operation research and computer science, such as agriculture, transportation, telecommunications, economics, production and operations scheduling, strategic decision-making, etc~\cite{hazell1974competitive,delson1992linear,dahleh1994control,liu2008choice,charnes1959application,zhou2008linear}. 

Since the 1940s, speeding up and scaling up LP has been a central topic in the optimization community, with extensive studies from both academia and industry. The current general-purpose LP solvers, such as Gurobi~\cite{optimization2023gurobi}, COPT~\cite{ge2022cardinal}, CPLEX~\cite{manual1987ibm} and HiGHS~\cite{huangfu2018parallelizing}, are quite mature. These classic LP solvers are based on either the simplex method or interior-point methods (IPMs), which can generally provide high-quality solutions to LP. However, further scaling up or speeding up these methods is highly challenging. The fundamental difficulty is the necessity of solving linear systems in simplex and IPMs, which requires either LU factorization (for simplex method) or Cholesky factorization (for IPMs). Such factorization-based approaches have two major drawbacks when solving large instances: (i) Storing the factorization can be quite memory-demanding. It is often the case that a sparse matrix has a much denser factorization, which is one of the reasons for classic solvers to raise an “out-of-memory” error even if the instance can be stored in memory; (ii) Both methods are highly nontrivial to exploit massive parallelization, {and thus they are challenging to take advantage of modern computing architectures such as graphic processing units (GPUs) and distributed computing, due to the sequential nature of factorizations}. Accordingly, all classic LP solvers are CPU-based and implemented on a single shared memory machine.

In contrast, the modern deep learning models used in practice, such as GPT-4, have trillions of variables, which are arguably significantly larger than the scale of LP commercial solvers are capable to solve\footnote{We comment that ``solving'' means a different order of accuracy tolerance in LP and deep learning models, and nevertheless, there is a significant gap in scale.}. 
One commonly believed reason for the success of deep learning models is the extensive use of GPUs and distributed computing, which massively speeds up the neural network training process. Indeed, commercial solver companies, such as Gurobi, always try to use GPUs to speed up their solvers, but the early efforts were unsuccessful~\cite{gurobigpu2015,gurobigpu2023}. The fundamental reason is that GPUs do not work well for solving sparse linear systems, which is the computational bottleneck of simplex or barrier method solving linear programming~\cite{swirydowicz2022linear,gurobigpu2015, gurobigpu2023}\footnote{NVIDIA recently released a library, cuDSS, for direct solving of sparse linear systems on GPUs in March 2024.}. 

This paper revisits this natural question:
\begin{center}
    \textit{Are GPUs useful in solving LP?}
\end{center}
We provide an affirmative answer to this question by presenting cuPDLP.jl, a GPU implementation of restarted primal-dual hybrid gradient (PDHG) implemented in Julia. We present an extensive numerical study comparing cuPDLP.jl with its CPU counterpart and three methods implemented in the highly-optimized commercial LP solver Gurobi on standard LP benchmark sets, which showcase that
\begin{itemize}
    \item cuPDLP.jl, a GPU implementation of PDLP in Julia programming language, has a comparable behavior with the highly-optimized commercial LP solver Gurobi on standard LP benchmark sets.
    \item Compared with its CPU counterpart, cuPDLP.jl clearly has superior performance on standard benchmark sets. Furthermore, one can observe the strong correlation between the GPU speed-up and the size of the instances.
\end{itemize}

cuPDLP.jl can be viewed as a CUDA implementation of its CPU counterpart PDLP~\cite{applegate2021practical}. PDLP has two open-sourced implementations, a prototype implementation in Julia (\href{https://github.com/google-research/FirstOrderLp.jl}{FirstOrderLp.jl}), and a production-level C++ implementation (open-sourced through \href{https://developers.google.com/optimization}{Google OR-Tools}).

Different from classic LP solvers, PDLP is a first-order method (FOM) LP solver. The fundamental difference is that the computational bottleneck of FOM is (sparse) matrix-vector multiplication, in contrast to (sparse) matrix factorization in simplex or IPMs. Table \ref{tab:compare} presents a comparison summary of simplex, barrier, and FOM LP solvers in five different dimensions: cost per iteration, number of iterations needed to solve an instance, code complexity, whether they can take advantage of massive parallelization, and whether they need extensive memory usage (beyond just storing the instance itself). As we can see, FOMs have multiple advantages in solving large instances, compared with simplex and IPMs. 

In particular, thanks to the recent development of deep learning, (sparse) matrix-vector multiplication suits very well on modern GPU infrastructure. Additionally, it is also shown that PDLP has strong theoretical guarantees~\cite{applegate2023faster,applegate2021infeasibility,lu2022infimal,lu2023geometry} and superior numerical performance compared to other FOM-based solvers~\cite{applegate2021practical}. These are the reasons we choose PDLP as the base algorithm in our GPU implementation.

To fully unleash the potential offered by modern GPU hardware, there are two major differences and modifications to the original CPU-based PDLP: (1) Due to the slow communication between CPU and GPU, the computational framework of PDLP has to be fully implemented on GPUs. In other words, both the instances and intermediate iterates must be resident within GPU memory (see Section \ref{sec:design} for details). (2) In terms of the algorithm, the only significant difference is a major restart scheme -- the CPU implementation utilizes normalized duality gap as the metric for restart~\cite{applegate2023faster}, while we utilize KKT error as the metric for restart (see Section \ref{sec:restart-kkt} for details). This avoids the trust-region algorithm to compute the normalized duality gap proposed in \cite{applegate2023faster}, which is sequential in nature and is not friendly for GPU implementation. Furthermore, we show in theory that restarting with KKT error enjoys the same order of linear convergence rate as restarting with normalized duality gap in Appendix \ref{app:thoery}.

This methodology is not limited to solving LP, and can potentially be extended to other optimization problems, such as quadratic programming~\cite{lu2023practical} and primal heuristics in mixed-integer programming~\cite{mexi2023scylla}.

\begin{table}
\centering
\includegraphics[width=0.7\textwidth]{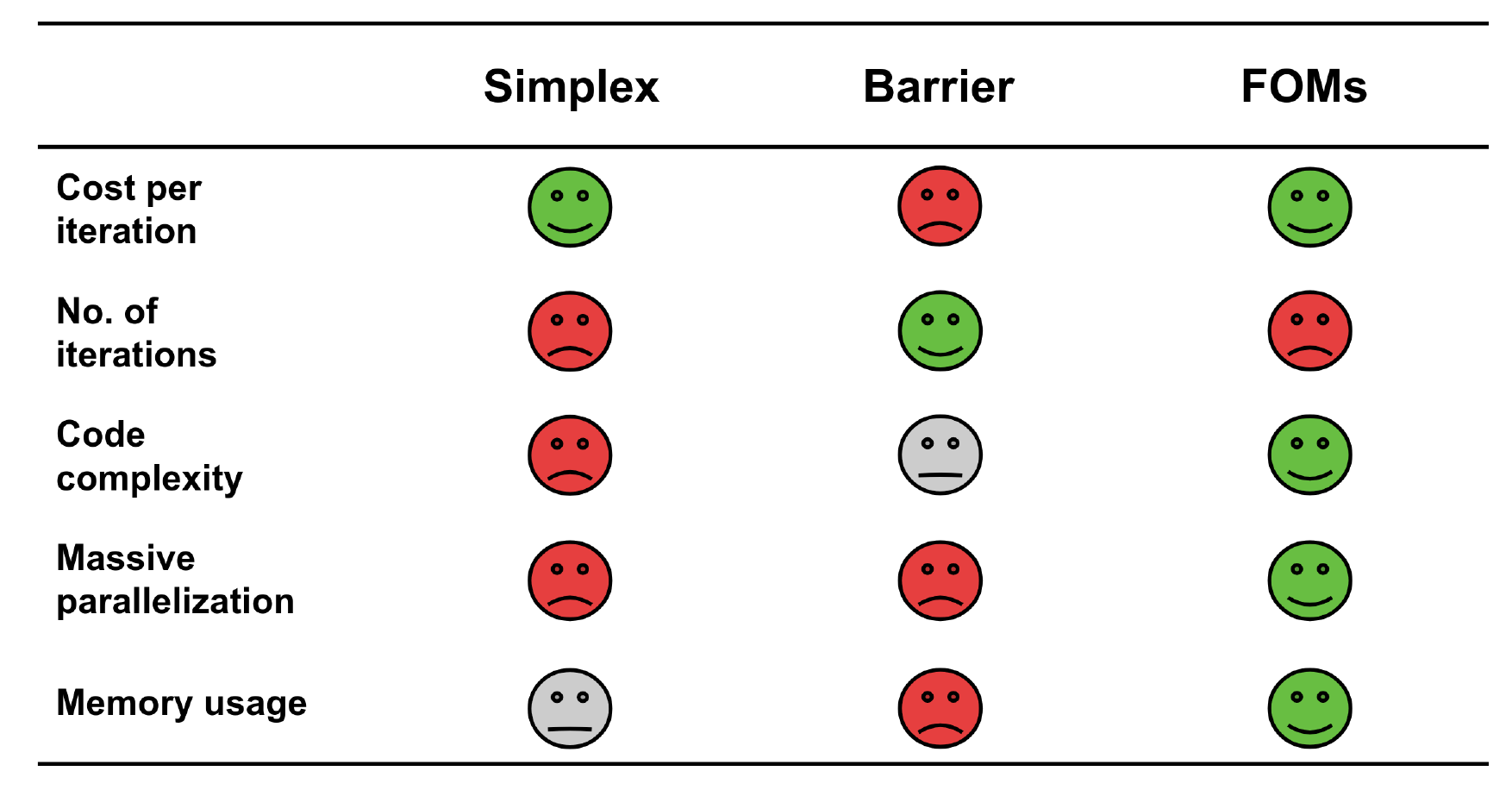}
\caption{Comparison of simplex, barrier, and first-order methods in five dimensions, where green, gray, and red faces represent a decreasing order of favorableness of the corresponding method in the corresponding dimension.} 
\label{tab:compare} 
\end{table}

\subsection{Related literature}
{\bf Linear programming and classic solvers.} 
The state-of-the-art methods to solve LP are simplex methods~\cite{dantzig1998linear} and interior-point methods \cite{karmarkar1984new}. Based on these methods, commercial solvers such as Gurobi~\cite{optimization2023gurobi} and COPT~\cite{ge2022cardinal} and open-sourced solvers like HiGHS~\cite{huangfu2018parallelizing} can provide reliable solutions with high accuracy.

{\bf FOM-based LP solvers.} There is a recent surge of research on using first-order methods (FOMs) to solve linear programming. FOMs are appealing due to their low per-iteration cost and ability of parallelization.
\begin{itemize}
    \item PDHG-based solvers: \href{https://github.com/google/or-tools/tree/stable/ortools/pdlp}{PDLP}~\cite{applegate2021practical}. PDLP is a general-purpose large-scale LP solver. PDLP is built upon restarted PDHG algorithm~\cite{applegate2023faster}, with many practical algorithmic enhancements, such as preconditioning, adaptive restart, adaptive step-size, infeasibility detection, etc. Currently PDLP has three implementations: a prototype implemented in Julia (\href{https://github.com/google-research/FirstOrderLp.jl}{FirstOrderLp.jl}), a production-level C++ implementation (open-sourced through \href{https://developers.google.com/optimization}{Google OR-Tools}), and an internal distributed version at Google~\cite{blog}. There is also an extensive study on the theoretical results related to PDLP, such as, the linear convergence rate~\cite{lu2022infimal,applegate2023faster}, infeasibility detection~\cite{applegate2021infeasibility}, refined complexity analysis~\cite{lu2023geometry,hinder2023worst}, extensions to quadratic programming~\cite{lu2023practical}, etc.
    \item Matrix-free IPM solvers: \href{https://github.com/leavesgrp/ABIP}{ABIP}~\cite{lin2021admm, deng2022new}. The core algorithm of ABIP is solving the homogeneous self-dual embedded cone programs via an interior-point method, and as a special case, it can solve LP. ABIP utilizes multiple ADMM iterations instead of one Newton step to approximately minimize the log-barrier penalty function. A recently enhanced version of ABIP (named ABIP+~\cite{deng2022new}) includes many new enhancements, such as preconditioning, restart, and hybrid parameter tuning, on top of ABIP. It was shown that ABIP+ (developed in C) has a comparable numerical performance on the LP benchmark sets to the Julia implementation of PDLP.

    \item Dual-based solvers: ECLIPSE~\cite{basu2020eclipse}. ECLIPSE is a distributed LP solver. It leverages accelerated gradient descent to solve a smoothed dual form of LP. ECLIPSE is designed specifically to solve large-scale LPs with certain decomposition structures arising from web applications. For example, ECLIPSE is used to solve real-world web applications with $10^{12}$ decision variables at LinkedIn~\cite{ramanath2022efficient,acharya2023promoting}.

    \item ADMM-based solvers: \href{https://github.com/cvxgrp/scs}{SCS}~\cite{o2016conic,o2021operator} and \href{https://github.com/osqp/osqp}{OSQP}~\cite{stellato2020osqp}. SCS is designed to solve convex cone programs, and OSQP is designed to solve convex quadratic programs. Both solvers are based on ADMM, and the major algorithmic difference is that SCS solves the homogeneous self-dual embedding of general conic programming, while OSQP solves the primal-dual form of quadratic program directly.
    LP can be solved as a special case of cone programming in SCS as well as a special case of quadratic programming in OSQP. The computational bottleneck of ADMM-based methods is solving a linear system with similar forms every iteration. 
    Both SCS and OSQP support using either direct solve of the linear system via factorization or indirect solve of the linear system via the conjugate gradient method. The first approach still requires doing one factorization and thus may suffer from the same scaling difficulties as simplex and IPMs. The second approach usually requires several conjugate gradient steps (i.e., matrix-vector multiplications) for every iteration.
\end{itemize}

{\bf Primal-dual hybrid gradient method (PDHG).} PDHG is an operator splitting method initially designed for applications in image processing~\cite{chambolle2011first,condat2013primal,esser2010general,he2012convergence,zhu2008efficient}. PDHG exhibits a sub-linear rate on general convex-concave problems~\cite{chambolle2016ergodic,lu2023unified} and achieves linear rate on a wide class of problems~\cite{lu2022infimal,fercoq2021quadratic}. PDLP~\cite{applegate2021practical} utilizes PDHG as its base algorithm. Built on the ergodic convergence of PDHG, a restart variant of PDHG exhibits an optimal linear convergence rate for solving LP~\cite{applegate2023faster}. Furthermore, Applegate et al. \cite{applegate2021infeasibility} shows how to extract infeasibility information of LP from PDHG iterates. 

{\bf Optimization with GPUs.} Previous efforts of commercial solver companies to solve LP on GPUs have not yet been successful due to the inefficiency of solving sparse linear systems on GPUs~\cite{gurobigpu2015,gurobigpu2023,swirydowicz2022linear}.   Recently, there have been open-sourced efforts using GPUs to solve generic optimization problems. For example, 
the ADMM-based solvers SCS~\cite{o2021operator} and OSQP~\cite{stellato2020osqp} have an implementation of indirect solves via conjugate gradient on GPUs, which avoids solving linear equations. \href{https://github.com/MadNLP/MadNLP.jl}{MadNLP.jl}~\cite{shin2023accelerating} is an interior-point method based solver for nonlinear programming that can be run on GPUs. However, these solvers are designed for solving different and more general classes of optimization problems, and they may not be suitable for solving LP instances with the size we consider herein (see Appendix \ref{app:comparison} for a comparison with SCS).

\section{PDLP}
In this section, we introduce the basics of LP and the algorithmic components of PDLP, as stated in~\cite{applegate2021practical}. After a brief introduction of vanilla PDHG for solving LP (Section \ref{sec:pdhg}), we summarize the enhancements of PDLP on top of PDHG that boost the practical performance (Section \ref{sec:enhancement}).

Consider LP of the general form:
\begin{equation}\label{eq:lp}
    \begin{aligned}[c]
    \min_{x\in \mathbb R^n}~~ &~ c^\top x \\
\text{s.t.}~~ &~ Gx \geq h \\
& ~ Ax = b \\
& ~ l \leq x \leq u\ ,
    \end{aligned}
\end{equation} 
where $G \in \mathbb R^{m_1\times n}$, $A \in \mathbb R^{m_2 \times n}$, $c \in \mathbb R^{n}$, $h \in \mathbb R^{m_1}$, $b \in \mathbb R^{m_2}$, $l \in (\mathbb R \cup \{ -\infty \})^{n}$, $u \in (\mathbb R \cup \{ \infty \})^{n}$. Notice that the constraint set is a complicated polytope that is hard to be projected onto, which makes it intractable to use standard FOMs, such as, projected gradient descent for solving \eqref{eq:lp}. PDLP solves the primal-dual form of the problem by dualizing the linear constraints:
\begin{flalign}\label{eq:primal-dual}
\min_{x \in X} \max_{y \in Y}\ L(x,y) := c^\top x - y^\top K x + q^\top y \ ,
\end{flalign}
where $K^\top = \begin{pmatrix} G^\top, A^\top \end{pmatrix}$ and $q^\top := \begin{pmatrix}
h^\top,
b^\top
\end{pmatrix}$, $X := \{x \in \mathbb R^n : l \leq x \leq u \}$, and $Y := \{y \in \mathbb R^{m_1+m_2} : y_{1:m_1} \geq 0\}.$ By duality theory, it is straightforward to show that a saddle point of \eqref{eq:primal-dual} recovers an optimal primal-dual solution to \eqref{eq:lp}.

\subsection{PDHG for solving LP}\label{sec:pdhg}
The base algorithm of PDLP is the primal-dual hybrid gradient (PDHG, a.k.a. Chambolle-Pock method)~\cite{chambolle2011first}. The update of PDHG for solving \eqref{eq:primal-dual} is 
\begin{equation}\label{eq:pdhg}
    \begin{cases}
        x^{t+1}\leftarrow \text{proj}_{X}(x^t-\tau (c-K^\top y^t)) \\ y^{t+1}\leftarrow \text{proj}_{Y}(y^t+\sigma (q-K(2x^{t+1}-x^t)))\ ,
    \end{cases}
\end{equation}
where $\tau,\sigma$ are primal and dual step-sizes respectively. The computational bottleneck of PDHG is matrix-vector multiplication $K^\top y^t$ and $Kx^t$; thus, PDHG is fully matrix-free to solve LP, i.e., it does not need to solve any linear equations. The primal and the dual step-size are reparameterized in PDLP as
\begin{equation*}
    \tau = \eta/\omega,\; \sigma=\eta\omega\quad \text{with}\;\  \eta,\omega>0\ ,
\end{equation*}
where $\eta$ (called step-size) controls the scale of the step-sizes, and $\omega$ (called primal weight) balances the primal and the dual progress. 

\subsection{Algorithmic enhancements in PDLP}\label{sec:enhancement}
We here summarize the algorithmic enhancements in PDLP that is presented in \cite{applegate2021practical}.
It turns out that the numerical performance of vanilla PDHG \eqref{eq:pdhg} on LP is not strong enough to support a modern solver~\cite{applegate2021practical}. To boost the practical performance, PDLP has essentially five major enhancements on top of PDHG: preconditioning, adaptive restart, adaptive step-size, primal weight update, and infeasibility detection. The core algorithm in PDLP, i.e., restarted PDHG, is presented in Algorithm \ref{alg:pdhg-restart}, followed by discussions on each algorithmic enhancement adopted.
\begin{algorithm}[h!]
\caption{Restarted PDHG (after preconditioning)}
\label{alg:pdhg-restart}
\SetKwInOut{Input}{Input}
\Input {Initial point $z^{0,0}$\;}
Initialize outer loop counter $n\leftarrow 0$, total iterations $k\leftarrow 0$, step-size $\hat\eta^{0,0}\leftarrow 1/\|K\|_{\infty}$, primal weight $\omega^0\leftarrow \text{InitializePrimalWeight}(c,q)$\;
\Repeat{\upshape termination criteria holds}{
  $t\leftarrow 0$\;
  \Repeat{\upshape restart or termination criteria holds}{
    $z^{n,t+1},\eta^{n,t+1},\hat\eta^{n,t+1}\leftarrow\text{AdaptiveStepPDHG}(z^{n,t},\omega^n,\hat\eta^{n,t},k)$\;
    $\bar z^{n,t+1}\leftarrow\tfrac{1}{\sum_{i=1}^{t+1}\eta^{n,i}}\sum_{i=1}^{k+1} \eta^{n,i}z^{n,i}$\;
    $z_c^{n,t+1}\leftarrow\text{GetRestartCandidate}(z^{n,t+1},\bar z^{n,t+1})$\;
    $t\leftarrow t+1$, $k\leftarrow k+1$\;
  }
  restart the outer loop $z^{n+1,0}\leftarrow z_c^{n,t}$, $n\leftarrow n+1$\;
  $\omega^n\leftarrow \text{PrimalWeightUpdate}(z^{n,0},z^{n-1,0},\omega^{n-1})$\;
}
\SetKwInOut{Output}{Output}
\Output{$z^{n,0}$.}
\end{algorithm}

\begin{itemize}
    \item {\bf Preconditioning:} The efficacy of first-order methods is closely tied to the conditioning of the underlying problem. In order to mitigate the ill-posedness, PDLP employs a diagonal preconditioner to ameliorate the condition number of the original problem. Specifically, it involves the rescaling of the constraint matrix $K=(G,A)$ to $\tilde K=(\tilde G,\tilde A)=D_1KD_2$ where $D_1$ and $D_2$ are positive diagonal matrices. This rescaling ensures that the resulting matrix $\tilde{K}$ is "well balanced". Consequently, this preconditioning step gives rise to a modified LP instance, wherein $A,G,c,b,h,u$ and $l$ in \eqref{eq:lp} are replaced with $\tilde G,\tilde A,\hat x=D_2^{-1}x,\tilde c=D_2c,(\tilde b,\tilde h)=D_1(b,h), \tilde u=D_2^{-1}u$ and $\tilde l=D_2^{-1}l$. In the default PDLP configuration, a combination of Ruiz rescaling \cite{ruiz2001scaling} and the preconditioning technique proposed by Pock and Chambolle \cite{pock2011diagonal} is employed.

    \item {\bf Adaptive restarts.} 
    PDLP utilizes an adaptive restarting strategy to enhance convergence. PDLP initially selects a restart candidate at each iteration, choosing between the current iterate and the average iterate based on a greedy principle. Subsequently, various restart criteria are assessed to determine if there is a constant factor decay in the progress metric. If such decay is observed, a restart is triggered. Further details can be found in \cite{applegate2023faster}.

    In the CPU-based PDLP, the progress metric for restarting is the normalized duality gap proposed in \cite{applegate2023faster}, and a trust-region algorithm is devised to compute this metric efficiently. While the trust-region algorithm exhibits linear time complexity, it is less compatible with GPUs' massively parallel computing paradigm. Therefore, in cuPDLP.jl, we introduce a novel restart scheme based on the KKT error. Detailed discussions are deferred to Section \ref{sec:restart-kkt}.
 
    \item {\bf Adaptive step-size.} The step-size suggested by theoretical considerations, namely $1/\|A\|_2$, turns out to be conservative in practical applications. To address this, PDLP employs a heuristic line search to determine a suitable step-size satisfying the condition:
    \begin{equation}\label{eq:step-size}
            \eta\leq \frac{\|z^{t+1}-z^t\|_{\omega}^2}{2(y^{t+1}-y^t)^\top K(x^{t+1}-x^t)}\ ,
    \end{equation}
    where $\|z\|_{\omega}:=\sqrt{\omega\|x\|_2^2+\frac{\|y\|_2^2}{\omega}}$ and $\omega$ is the current primal weight. Additional details of the adaptive step-size rule are elaborated in \cite{applegate2021practical}. The inequality \eqref{eq:step-size} was inspired from the $\mathcal O(1/k)$ convergence rate proof of PDHG~\cite{chambolle2016ergodic,lu2023unified}. The empirical evidence from numerical experiments conducted in \cite{applegate2021practical} attests to its consistent efficacy.
    \begin{algorithm}
    \caption{One step of PDHG using adaptive step-size heuristic}
    \label{alg:adaptive-step}
    \SetKwInOut{Input}{Function}
    \Input{AdaptiveStepPDHG($z^{n,t}$, $\omega^n$, $\hat\eta^{n,t}$,$k$)}

        $(x,y)\leftarrow z^{n,t}$, $\eta\leftarrow \hat\eta^{n,t}$\\
        \For{$i=0,1,...$}{
        $x'\leftarrow \mathrm{proj}_X(x-\frac{\eta}{\omega^n}(c-K^\top y))$\\
        $y'\leftarrow \mathrm{proj}_Y(y+{\eta}{\omega^n}(q-K(2x'-x)))$\\
        $\bar \eta\leftarrow\frac{\|(x'-x,y'-y)\|_{\omega^n}^2}{2(y'-y)^\top K(x'-x)}$\\
        $\eta'\leftarrow\min\pran{(1-(k+1)^{-0.3})\bar \eta,(1+(k+1)^{-0.6})\eta}$\\
        \If{$\eta\leq \bar \eta$}{
        return $(x',y')$, $\eta$, $\eta'$
        }
        $\eta\leftarrow \eta'$
        }     
\end{algorithm}

\item {\bf Primal Weight Update:}
Adjusting the primal weight $\omega$ is designed to harmonize the primal and dual spaces through a heuristic approach. The update of primal weight is specific during restart occurrences, thus infrequently. More precisely, the initialization of $\omega$ involves the expression:
\begin{equation*}
\mathrm{InitializePrimalWeight}(c,q):=\begin{cases}
    \frac{\|c\|_2}{\|q\|_2},\; & \text{if } \|c\|_2,\|q\|_2>\epsilon_{\mathrm{zero}}\\
    1, \; & \mathrm{otherwise}
\end{cases}
\end{equation*}
where $\epsilon_{\mathrm{zero}}$ denotes a small nonzero tolerance. Let $\Delta_x^n=\|x^{n,0}-x^{n-1,0}\|_2$ and $\Delta_y^n=\|y^{n,0}-y^{n-1,0}\|_2$. PDLP initiates the primal weight update at the beginning of each new epoch.
{\small
\begin{equation*}
\mathrm{PrimalWeightUpdate}(z^{n,0},z^{n-1,0},\omega^{n-1}):=\begin{cases}
    \exp\pran{\theta \log\pran{\frac{\Delta_y^n}{\Delta_x^n}}+(1-\theta)\omega^{n-1}},\; & \Delta_x^n,\Delta_y^n>\epsilon_{\mathrm{zero}}\\
    \omega^{n-1}, \; & \mathrm{otherwise}
\end{cases}
\end{equation*}
}

The intuition is to determine the primal weight $\omega^n$ in a manner that equalizes the distance to optimality in both the primal and dual domains, i.e., $\|(x^{n,t}-x^*,0)\|_{\omega^n}\approx \|(0,y^{n,t}-y^*)\|_{\omega^n}$. Additionally, PDLP employs exponential smoothing with a parameter $\theta \in [0,1]$ to mitigate oscillations.
    
    \item {\bf Infeasibility detection.} PDLP periodically checks whether the difference of iterates $z^{n,t+1}-z^{n,t}$ or the normalized iterates $\frac{1}{t}(z^{n,t}-z^{n,0})$ provide an infeasibility certificate, and the performance of these two sequences is instance-dependent. A detailed investigation of infeasibility detection using PDHG iterates is available in \cite{applegate2021infeasibility}.
\end{itemize}

\section{GPU implementation of PDLP}
This section presents the design of cuPDLP.jl, the GPU implementation of PDLP. Section \ref{sec:architecture} briefly introduces the hardware architecture and logical structure. In Section \ref{sec:design}, we discuss in detail the design of cuPDLP.jl. Section \ref{sec:operations} presents the implementation of basic matrix and vector operations in cuPDLP.jl, and Section \ref{sec:restart-kkt} shows the restart scheme based on KKT error to adapt on GPU. The solver is available at \href{https://github.com/jinwen-yang/cuPDLP.jl}{https://github.com/jinwen-yang/cuPDLP.jl}.

\subsection{GPU architecture and thread hierarchy}\label{sec:architecture}
GPUs exhibit distinct underlying architecture with CPUs.  GPUs have significantly more computation cores than CPUs. For example, the GPU we used in the experiments, NVIDIA H100, has 7296 double-precision cores. However, unlike CPU cores, many GPU cores share the same control unit and must execute the same instruction simultaneously. As a result, the hardware design of GPUs encourages high bandwidth instead of a deep pipeline.

GPUs follow the single instruction multiple data (SIMD) computational paradigm; namely, the threads execute the same instruction but fetch their own data. At the heart of GPU programming is the kernel function, i.e., a program designed to execute instructions on GPUs. Upon launching the kernel function, the execution environment configures a grid of thread blocks, each block comprising an identical number of threads. A block of threads is assigned to an available streaming multiprocessor (SM) during runtime, directing them into warps, with each warp typically encompassing a set of 32 threads. Warp is the basic unit of execution in GPUs. One caveat goes that this does not mean that all thread blocks can run concurrently on SMs, and there are no guarantees on the order of block and warp execution.

\begin{figure}[ht!]
 \centering
     \includegraphics[width=0.5\textwidth]{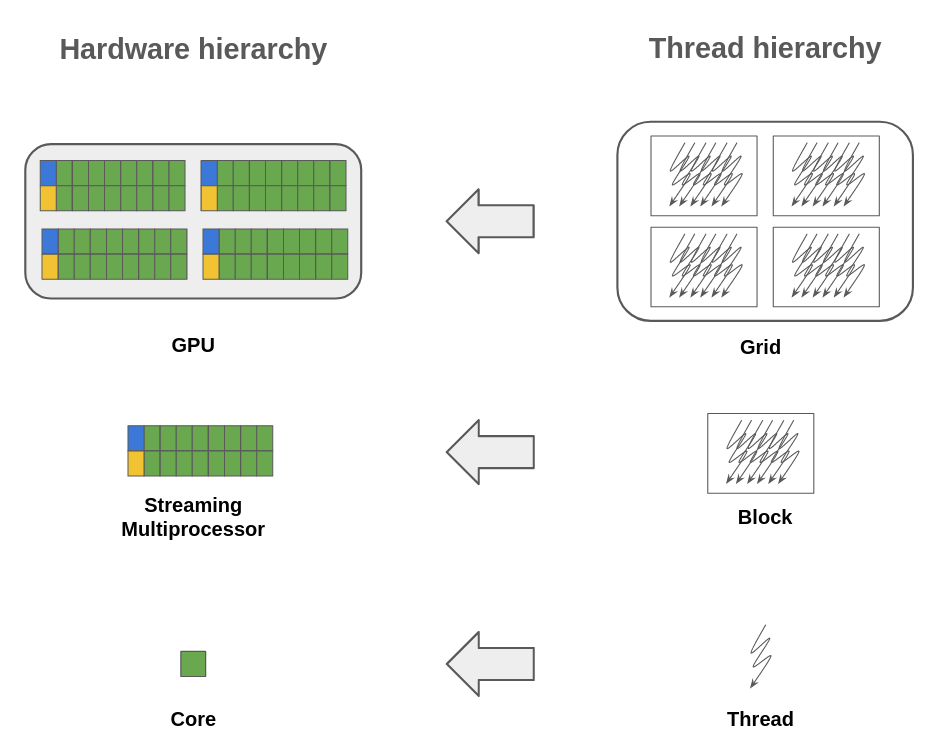}
	\caption{Illustration of the relationship between physical architecture and logical structure.}
	\label{fig:hierarchy}
\end{figure}

This hierarchical structure is the backbone of GPUs to foster parallel execution across threads, {which improves computational efficiency by reducing the latency and increasing the throughput (the data processed per time unit).} The relationship between the physical architecture and logical structure is depicted in Figure \ref{fig:hierarchy}. For more detailed discussion in the characteristic of GPUs, refer to \cite{guide2023cuda}.

Each GPU device has its own memory hardware, which is separate from the memory of the CPU host. A typical computational paradigm is using CPUs for IO process and exploiting GPUs for burdensome computation. Thus, data moving between CPU and GPU are needed. However, the significant cost of CPU-GPU communication, especially for large-scale problems, can be a crucial issue. Frequent data transfers between CPU and GPU can introduce dominating communication cost over any computational speed-up gained on GPU.

\subsection{Design of cuPDLP.jl}\label{sec:design}
\begin{figure}[ht!]
	\centering
    \includegraphics[width=0.5\textwidth]{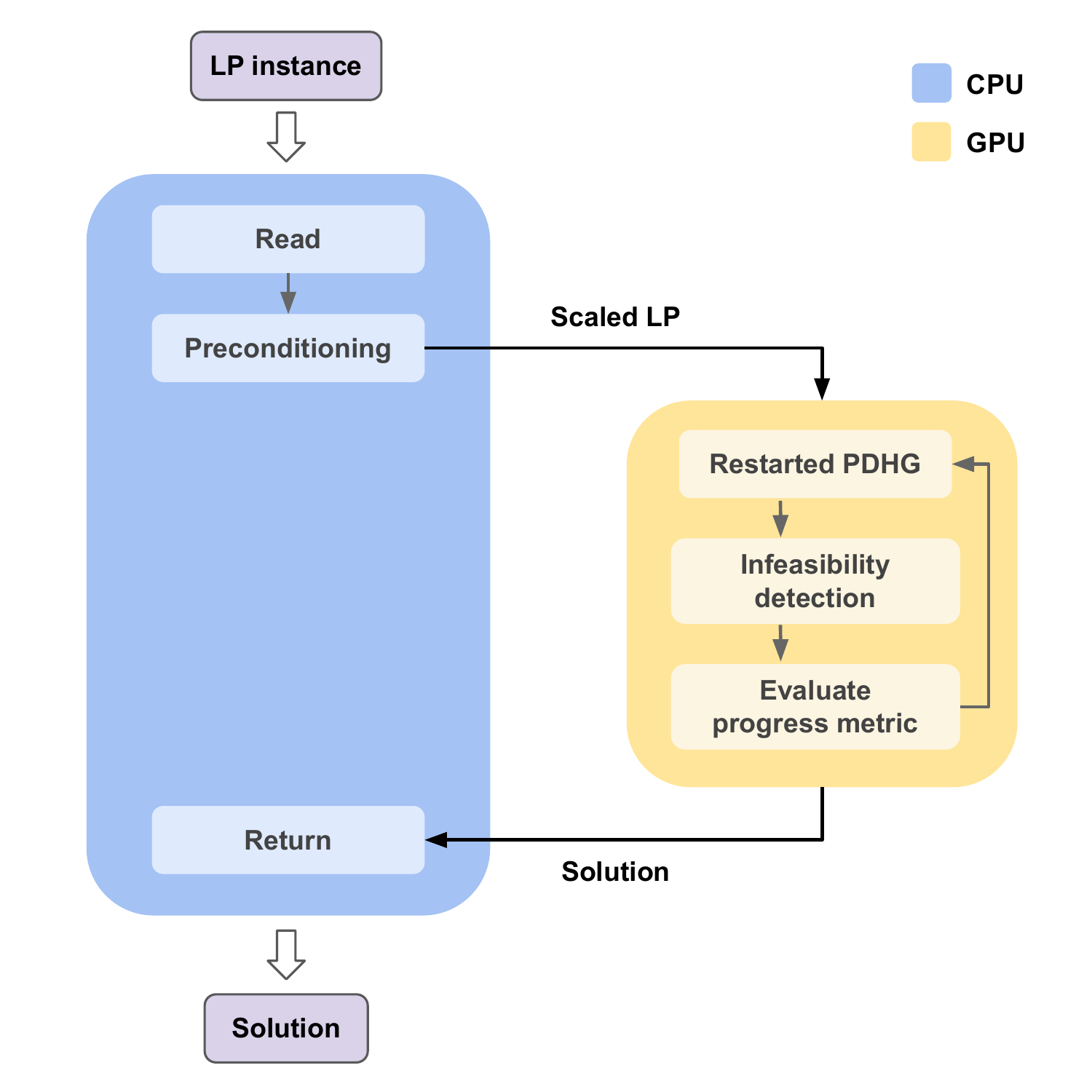}
	\caption{Illustration of the computational architecture of cuPDLP.jl.}
	\label{fig:diagram}
\end{figure}

The design of cuPDLP.jl is illustrated in Figure \ref{fig:diagram}. To avoid expensive data transfers between CPU and GPU, we have designed our implementation of cuPDLP.jl to run as much as possible on the GPU. 

As depicted in Figure \ref{fig:diagram}, only two communications between CPU and GPU are required for cuPDLP.jl. One transfers the scaled LP instance after preconditioning from CPU to GPU while the other moves the final solution from GPU to CPU as output. Major iterations are executed completely on GPU and there is no need to transfer any vectors before termination to alleviate expensive CPU-GPU communication.

\subsection{Vector and matrix operations}\label{sec:operations}
Matrix-vector multiplications and vector-vector operations, the core of PDLP, can be parallelized naturally and fit well with the SIMD paradigm of GPUs. GPUs have massive parallelization capability to manipulate each coordinate of vectors on each of its core in a parallel fashion. Specifically, the constraint matrix is stored in Compressed Sparse Row (CSR) format and cuPDLP.jl uses the 1-dimensional thread configuration since most of our operations are matrix-vector
and vector-vector operations. We write our custom kernels for main PDHG updates \eqref{eq:pdhg} in cuPDLP.jl and utilize cusparseSpMV() implemented in cuSPARSE library to do matrix-vector multiplications. In our customized kernel, each thread updates a single coordinate of the iterate vector to maximize the throughput and cusparseSpMV() uses algorithm \texttt{CUSPARSE\_SPMV\_CSR\_ALG2} to provide a deterministic result for each run.

\subsection{Adaptive restart based on KKT error}\label{sec:restart-kkt}
PDLP utilizes normalized duality gap as the progress metric for restart, which was introduced in~\cite{applegate2023faster}. The normalized duality gap for LP is computed by a trust-region algorithm. Although only linear time is needed, it is nontrivial to implement the trust-region method in an efficient parallel version on GPU due to the sequential nature of the trust-region method. We use the KKT error defined in \eqref{eq:kkt} as a proxy to the normalized duality gap for restarting since KKT error can be evaluated very efficiently on GPU.
\begin{equation}\label{eq:kkt}
    \mathrm{KKT}_{\omega}(z)=\sqrt{\omega^2\left\|\begin{pmatrix}
        Ax-b \\ [h-Gx]^+
    \end{pmatrix}\right\|_2^2+\frac{1}{\omega^2}\|c-K^\top y-\lambda\|_2^2+(q^\top y+l^\top\lambda^+-u^\top\lambda^--c^\top x)^2}
\end{equation}
More specifically, the restart scheme is described as follows

{\bf Choosing the restart candidate.} 
\begin{equation*}
    z_c^{n,t+1}:=\mathrm{GetRestartCandidate}(z^{n,t+1},\bar z^{n,t+1})=\begin{cases}
        z^{n,t+1},& \mathrm{KKT}_{\omega^n}(z^{n,t+1})<\mathrm{KKT}_{\omega^n}(\bar z^{n,t+1}) \\ \bar z^{n,t+1},& \mathrm{otherwise} \ .
    \end{cases}
\end{equation*}

{\bf Restart criteria.} Define parameters $\beta_{\mathrm{sufficient}}=0.2$, $\beta_{\mathrm{necessary}}=0.8$ and $\beta_{\mathrm{artificial}}=0.36$. Denote $k$ the total iteration counter. The algorithm restarts if one of three conditions holds:
\begin{enumerate}
    \item[(i)] (Sufficient decay in KKT error)
    \begin{equation*}
        \mathrm{KKT}_{\omega^n}(z_c^{n,t+1})\leq \beta_{\mathrm{sufficient}}\mathrm{KKT}_{\omega^n}(z^{n,0}) \ ,
    \end{equation*}
    \item[(ii)] (Necessary decay + no local progress in KKT error)
    \begin{equation*}
        \mathrm{KKT}_{\omega^n}(z_c^{n,t+1})\leq \beta_{\mathrm{necessary}}\mathrm{KKT}_{\omega^n}(z^{n,0})\quad \text{and} \quad \mathrm{KKT}_{\omega^n}(z_c^{n,t+1})>\mathrm{KKT}_{\omega^n}(z_c^{n,t}) \ ,
    \end{equation*}
    \item[(iii)] (Long inner loop)
    \begin{equation*}
        t\geq \beta_{\mathrm{artificial}}k \ .
    \end{equation*}
\end{enumerate}
Finally, we highlight that the only difference between the restart scheme here and that in \cite{applegate2021practical,applegate2023faster} is a different choice of metric: we utilize KKT residual, while \cite{applegate2021practical,applegate2023faster} utilized normalized duality gap.

\section{Numerical experiments}
In this section, we study the numerical performance of cuPDLP.jl. In particular, we compare cuPDLP.jl with the CPU implementations of PDLP~\cite{applegate2021practical} and three methods implemented in the commercial solver Gurobi~\cite{optimization2023gurobi}, i.e., primal simplex method, dual simplex method, and interior-point method. Section \ref{sec:setup} describes the setup of the experiments. Section \ref{sec:miplib} presents the numerical results on LP relaxations of instances from MIPLIB 2017 collection~\cite{gleixner2021miplib}. More specifically, Section \ref{sec:miplib-gurobi} investigates the performance of cuPDLP.jl and Gurobi, and Section \ref{sec:miplib-pdlp} compares cuPDLP.jl with different versions of PDLP, namely Julia implemented FirstOrderLp.jl and C++ implemented PDLP with single thread and multiple threads.  Section \ref{sec:mittelmann} discusses the performances of different solvers on Mittelmann's LP benchmark set~\cite{mittelmann}. We present additional figures and discussions on the results in Appendix \ref{app:figure}.

\subsection{Experimental setup}\label{sec:setup}

{\bf Benchmark datasets.} We use two LP benchmark datasets in the numerical experiments, \texttt{MIP Relaxations}, which contain 383 instances curated from root-node LP relaxation of mixed-integer programming problems from MIPLIB 2017 collection~\cite{gleixner2021miplib} (see Section \ref{sec:miplib}), and 49 LP instances from the \texttt{Mittelmann's LP} benchmark dataset~\cite{mittelmann} (see Section \ref{sec:mittelmann}).

In particular, MIPLIB 2017~\cite{gleixner2021miplib} is a collection of mixed-integer linear programming problems. We utilize the root-node LP relaxation of instances in MIPLIB as the LP benchmark set. 383 instances are selected from MIPLIB 2017 to construct \texttt{MIP Relaxations} based on the following criteria (similar selection criteria are used in the experiments of CPU-based PDLP~\cite{applegate2021practical}):
\begin{itemize}
    \item Not tagged as numerically unstable
    \item Not tagged as infeasible
    \item Not tagged as having indicator constraints
    \item Finite optimal objective (if known)
    \item The constraint matrix has a number of nonzeros greater than 100,000
    \item Zero is not an optimal solution to the LP relaxation.
\end{itemize}
\texttt{MIP Relaxations} is further split into three classes based on the number of nonzeros (nnz) in the constraint matrix, as shown in Table \ref{tab:miplib-size}.
\begin{table}[ht!]
\centering
\begin{tabular}{cccc}
\hline
                    & \textbf{Small}                   & \textbf{Medium}                     & \textbf{Large}                   \\ \hline
\textbf{Number of nonzeros}  & 100K -  1M & 1M - 10M & \textgreater 10M \\
\textbf{Number of instances} & 269                     & 94                         & 20                      \\ \hline
\end{tabular}
\caption{Scales of instances in \texttt{MIP Relaxations}.}
\label{tab:miplib-size}
\end{table}

\texttt{Mittelmann's LP} is a classic dataset to benchmark LP solvers. We utilize the 49 public instances from the dataset. The number of nonzeros of the instances spans from 100 thousand to 90 million.

{\bf Software.} \href{https://github.com/jinwen-yang/cuPDLP.jl}{cuPDLP.jl} is implemented in an open-source Julia~\cite{bezanson2017julia} module. cuPDLP.jl utilizes \href{https://github.com/JuliaGPU/CUDA.jl}{CUDA.jl}~\cite{besard2018effective} as the interface for working with NVIDIA CUDA GPUs using Julia. We compare cuPDLP.jl with five LP solvers: the Julia implementation of PDLP (\href{https://github.com/google-research/FirstOrderLp.jl}{FirstOrderLp.jl}), the C++ implementation of PDLP with single thread and multiple threads (wrapped in \href{https://github.com/google/or-tools}{Google OR-Tools}\footnote{Specifically, PDLP in \href{https://github.com/google/or-tools/tree/master}{master branch of Google OR-Tools} is used in our experiments.}), the primal simplex, dual simplex and barrier method implemented in Gurobi. Crossover is disabled, and 16 threads are used in experiments of Gurobi. {The running time of cuPDLP.jl and FirstOrderLp.jl is measured after pre-compilation in Julia.}

{\bf Computing environment.} We use NVIDIA H100-PCIe-80GB GPU, with CUDA 12.3, for running cuPDLP.jl, and we use Intel Xeon Gold 6248R CPU 3.00GHz with 160GB RAM and 16 threads for running CPU-based solvers.  The experiments are performed in Julia 1.9.2 and Gurobi 11.0 (released in November 2023). Table \ref{tab:flop} compares the theoretical peak of double-precision FLOPS (floating-point operations per second) of the CPU and GPU used in the experiments. 

\begin{table}[ht!]
\centering
{\small
\begin{tabular}{ccc}
\hline
                                      & \textbf{CPU (16 threads)}      & \textbf{GPU}                \\ \hline
\textbf{Processor}                    & Intel Xeon Gold 6248R CPU 3.00GHz\tablefootnote{See \href{https://ark.intel.com/content/www/us/en/ark/products/199351/intel-xeon-gold-6248r-processor-35-75m-cache-3-00-ghz.html}{processor specifications}.} & NVIDIA H100-PCIe\tablefootnote{H100 has 114 SMs and 7296 FP64 cores. See \href{https://resources.nvidia.com/en-us-tensor-core/nvidia-tensor-core-gpu-datasheet}{H100 datasheet} and \href{https://resources.nvidia.com/en-us-tensor-core/gtc22-whitepaper-hopper}{H100 whitepaper} for more a detailed description of H100 GPU. } \\
\textbf{Theoretical peak (FP64)} & 256 GFLOPS                         & 26 TFLOPS \\ 
\textbf{Maximum memory bandwidth} & 137.48 GB/sec  & 2 TB/sec                 \\ \hline
\end{tabular}
}
\caption{Comparison of CPU and GPU specifications.}
\label{tab:flop}
\end{table}

{\bf Initialization.} Both PDLP and cuPDLP.jl uses all-zero vectors as the initial starting points.

{\bf Optimality termination criteria.} PDLP and cuPDLP.jl terminate when the relative KKT error is no greater than the termination tolerance $\epsilon\in(0,\infty)$:
\begin{equation*}
    \begin{aligned}
        |q^\top y+l^\top \lambda^+-u^\top \lambda^--c^\top x|&\leq \epsilon(1+|q^\top y+l^\top\lambda^+-u^\top \lambda^-|+|c^\top x|)\\
        \left\|\begin{pmatrix}
            Ax-b \\ [h-Gx]^+
        \end{pmatrix}\right\|_2&\leq \epsilon(1+\|q\|_2)\\
        \|c-K^\top y-\lambda\|_2&\leq \epsilon(1+\|c\|_2) \ .
    \end{aligned}
\end{equation*}
The termination criteria are checked for the original LP instance, not the preconditioned ones, so that the termination is not impacted by the preconditioning.
We use $\epsilon=10^{-4}$ for moderately accurate solutions and $\epsilon=10^{-8}$ for high-quality solutions. We also set $10^{-4}$ and $10^{-8}$ tolerances for parameters \texttt{FeasibilityTol}, \texttt{OptimalityTol} and \texttt{BarConvTol} (for barrier methods) of Gurobi\footnote{We comment that Gurobi barrier, Gurobi simplex and PDLP all use different termination criteria, so it can never be a fully ``fair'' comparison among different types of algorithms.}.

{\bf Time limit.} In Section \ref{sec:miplib}, we impose a time limit of 3600 seconds on instances with small-sized and medium-sized instances and a time limit of 18000 seconds for large instances. For Section \ref{sec:mittelmann}, we impose 15000 seconds as the time limit in Mittelmann's benchmark~\cite{mittelmann}. 

{\bf Shifted geometric mean.} We report the shifted geometric mean of solve time to measure the performance of solvers on a certain collection of problems. More precisely, shifted geometric mean is defined as $\left(\prod_{i=1}^n (t_i+\Delta)\right)^{1/n}-\Delta$ where $t_i$ is the solve time for the $i$-th instance. We shift by $\Delta=10$ and denote it SGM10. If the instance is unsolved, the solve time is always set to the corresponding time limit. 

{\bf Presolves.} Presolve is a step that can simplify the LP problem before solving it. It involves detecting inconsistent bounds, removing empty rows and columns from the constraint matrix, eliminating variables with equal lower and upper bounds, detecting duplicate rows, tightening bounds, etc. This step is independent of the algorithm solving the instances. In the next section, we present the comparison with Gurobi on two sets of instances, namely, the original instances without any presolve step, as well as the instances after Gurobi presolve step.

\subsection{MIP Relaxations}\label{sec:miplib}
In this section, we compare cuPDLP.jl with the commercial LP solver Gurobi and CPU implementations of PDLP on \texttt{MIP Relaxations} respectively. The main message is that GPU-implemented cuPDLP.jl exhibits significant speedup over its CPU-implemented counterparts, and its numerical performance is on par with Gurobi.

\subsubsection{cuPDLP.jl versus Gurobi}\label{sec:miplib-gurobi}
\begin{table}[h!]
\centering
{\small
\begin{tabular}{ccccccccc}
\hline
\multirow{2}{*}{}                                   & \multicolumn{2}{c}{\begin{tabular}[c]{@{}c@{}}\textbf{Small (269)} \\ (1-hour limit)\end{tabular}} & \multicolumn{2}{c}{\begin{tabular}[c]{@{}c@{}}\textbf{Medium (94)}\\ (1-hour limit)\end{tabular}}    & \multicolumn{2}{c}{\begin{tabular}[c]{@{}c@{}}\textbf{\textbf{Large (20)}}\\ (5-hour limit)\end{tabular}}  & \multicolumn{2}{c}{\textbf{Total (383)}}       \\
                                                    & \textbf{Count} & \textbf{Time} & \textbf{Count} & \textbf{Time} & \textbf{Count} & \textbf{Time} & \textbf{Count} & \textbf{Time} \\ \hline
\multicolumn{1}{c}{\textbf{cuPDLP.jl}}     & 266                   & 8.61               & 92                    & 14.80               & 19                    & 111.19 &377 &12.02             \\
\multicolumn{1}{c}{\textbf{Primal simplex (Gurobi)}} & 268                   & 12.56                & 69                    & 188.81              & 11                    & 3145.49   &348 &39.81           \\
\multicolumn{1}{c}{\textbf{Dual simplex (Gurobi)}}   & 268                   & 8.75                & 84                    & 66.67               & 15                    & 591.63    &367 &21.75           \\
\multicolumn{1}{c}{\textbf{Barrier (Gurobi)}}        & 268                   & 5.30                & 88                    & 45.01               & 18                    & 415.78    &374 &14.92           \\ \hline
\end{tabular}
}
\caption{Solve time in seconds and SGM10 of different solvers on instances of \texttt{MIP Relaxations} with tolerance $10^{-4}$: cuPDLP.jl versus Gurobi without presolve.}
\label{tab:miplib-1e-4-no-presolve}
\end{table}

\begin{table}[h!]
\centering
{\small
\begin{tabular}{ccccccccc}
\hline
\multirow{2}{*}{}                                   & \multicolumn{2}{c}{\begin{tabular}[c]{@{}c@{}}\textbf{Small (269)} \\ (1-hour limit)\end{tabular}} & \multicolumn{2}{c}{\begin{tabular}[c]{@{}c@{}}\textbf{Medium (94)}\\ (1-hour limit)\end{tabular}}    & \multicolumn{2}{c}{\begin{tabular}[c]{@{}c@{}}\textbf{\textbf{Large (20)}}\\ (5-hour limit)\end{tabular}}  & \multicolumn{2}{c}{\textbf{Total (383)}}       \\
                                                    & \textbf{Count} & \textbf{Time} & \textbf{Count} & \textbf{Time} & \textbf{Count} & \textbf{Time} & \textbf{Count} & \textbf{Time} \\ \hline
\multicolumn{1}{c}{\textbf{cuPDLP.jl}}     & 269                   & 5.35               & 93                    & 10.31               & 19                    & 33.93 &381 &7.37             \\
\multicolumn{1}{c}{\textbf{Primal simplex (Gurobi)}} & 269                   & 5.67                & 71                    & 121.23              & 19                    & 297.59   &359 &20.84          \\
\multicolumn{1}{c}{\textbf{Dual simplex (Gurobi)}}   & 268                   & 4.17                & 86                    & 37.56               & 19                    & 179.49    &373 &11.84           \\
\multicolumn{1}{c}{\textbf{Barrier (Gurobi)}}        & 269                   & 1.21                & 94                    & 15.32               & 20                    & 30.70    &383 &4.65           \\ \hline
\end{tabular}
}
\caption{Solve time in seconds and SGM10 of different solvers on instances of \texttt{MIP Relaxations} with tolerance $10^{-4}$: cuPDLP.jl versus Gurobi with presolve.}
\label{tab:miplib-1e-4-with-presolve}
\end{table}

\begin{table}[h!]
\centering
{\small
\begin{tabular}{ccccccccc}
\hline
\multirow{2}{*}{}                                   & \multicolumn{2}{c}{\begin{tabular}[c]{@{}c@{}}\textbf{Small (269)} \\ (1-hour limit)\end{tabular}} & \multicolumn{2}{c}{\begin{tabular}[c]{@{}c@{}}\textbf{Medium (94)}\\ (1-hour limit)\end{tabular}}    & \multicolumn{2}{c}{\begin{tabular}[c]{@{}c@{}}\textbf{\textbf{Large (20)}}\\ (5-hour limit)\end{tabular}} & \multicolumn{2}{c}{\textbf{Total (383)}}          \\
                                                    & \textbf{Count} & \textbf{Time}  & \textbf{Count} & \textbf{Time}  & \textbf{Count} & \textbf{Time} & \textbf{Count} & \textbf{Time}  \\ \hline
\multicolumn{1}{c}{\textbf{cuPDLP.jl}}     & 261                   & 23.47                & 86                    & 40.69                & 16                    & 421.40  &363 &32.35             \\
\multicolumn{1}{c}{\textbf{Primal simplex (Gurobi)}} & 268                   & 12.43                 & 74                    & 157.59               & 13                    & 2180.23        &355 &36.68      \\
\multicolumn{1}{c}{\textbf{Dual simplex (Gurobi)}}   & 268                   & 8.00                 & 83                    & 59.93                & 15                    & 687.17       &366 &20.40      \\
\multicolumn{1}{c}{\textbf{Barrier (Gurobi)}}        & 267                   & 6.24                 & 88                    & 48.62                & 18                    & 438.69      &373 &16.46        \\ \hline
\multicolumn{1}{l}{}                                & \multicolumn{1}{l}{}  & \multicolumn{1}{l}{} & \multicolumn{1}{l}{}  & \multicolumn{1}{l}{} & \multicolumn{1}{l}{}  & \multicolumn{1}{l}{}
\end{tabular}
}
\caption{Solve time in seconds and SGM10 of different solvers on instances of \texttt{MIP Relaxations} with tolerance $10^{-8}$: cuPDLP.jl versus Gurobi without presolve.\tablefootnote{Gurobi primal and dual simplex methods may perform better for obtaining high-accuracy solution than medium-accuracy solution. This is a \href{https://support.gurobi.com/hc/en-us/community/posts/21790169951249-Simplex-method-does-not-stop-at-optimal-solution}{known effect} due to the different trajectory paths when setting different tolerance levels.}}
\label{tab:miplib-1e-8-no-presolve}
\end{table}

\begin{table}[h!]
\centering
{\small
\begin{tabular}{ccccccccc}
\hline
\multirow{2}{*}{}                                   & \multicolumn{2}{c}{\begin{tabular}[c]{@{}c@{}}\textbf{Small (269)} \\ (1-hour limit)\end{tabular}} & \multicolumn{2}{c}{\begin{tabular}[c]{@{}c@{}}\textbf{Medium (94)}\\ (1-hour limit)\end{tabular}}    & \multicolumn{2}{c}{\begin{tabular}[c]{@{}c@{}}\textbf{\textbf{Large (20)}}\\ (5-hour limit)\end{tabular}} & \multicolumn{2}{c}{\textbf{Total (383)}}          \\
                                                    & \textbf{Count} & \textbf{Time}  & \textbf{Count} & \textbf{Time}  & \textbf{Count} & \textbf{Time} & \textbf{Count} & \textbf{Time}  \\ \hline
\multicolumn{1}{c}{\textbf{cuPDLP.jl}}     & 264                   & 17.53                & 90                    & 30.05                & 19                    & 81.07  &373 &22.13             \\
\multicolumn{1}{c}{\textbf{Primal simplex (Gurobi)}} & 269                   & 5.19                & 75                    & 100.03               & 18                    & 171.72        &362 &18.11      \\
\multicolumn{1}{c}{\textbf{Dual simplex (Gurobi)}}   & 268                   & 3.53                 & 89                    & 27.17               & 19                    & 121.94        &376 &9.53      \\
\multicolumn{1}{c}{\textbf{Barrier (Gurobi)}}        & 269                   & 1.34                 & 94                    & 16.85                & 20                    & 33.48       &383 &5.03        \\ \hline
\multicolumn{1}{l}{}                                & \multicolumn{1}{l}{}  & \multicolumn{1}{l}{} & \multicolumn{1}{l}{}  & \multicolumn{1}{l}{} & \multicolumn{1}{l}{}  & \multicolumn{1}{l}{}
\end{tabular}
}
\caption{Solve time in seconds and SGM10 of different solvers on instances of \texttt{MIP Relaxations} with tolerance $10^{-8}$: cuPDLP.jl versus Gurobi with presolve.}
\label{tab:miplib-1e-8-with-presolve}
\end{table}

Table \ref{tab:miplib-1e-4-no-presolve}-\ref{tab:miplib-1e-8-with-presolve} present a comparison between cuPDLP.jl and the commercial LP solver Gurobi. Particularly, Table \ref{tab:miplib-1e-4-no-presolve} and Table \ref{tab:miplib-1e-4-with-presolve} present moderate accuracy results (i.e., $10^{-4}$ relative KKT error), while Table \ref{tab:miplib-1e-8-no-presolve} and Table \ref{tab:miplib-1e-8-with-presolve} present the high accuracy results (i.e., $10^{-8}$ relative KKT error); Table \ref{tab:miplib-1e-4-no-presolve} and Table \ref{tab:miplib-1e-8-no-presolve} present the results on solving the original problems, while Table \ref{tab:miplib-1e-4-with-presolve} and Table \ref{tab:miplib-1e-8-with-presolve} present the results on LP instances after Gurobi presolve. The tables yield several noteworthy observations:
\begin{itemize}
    \item With Gurobi presolve, cuPDLP.jl can solve 99.5\% instances to medium accuracy and 97.4\% instances to high accuracy within the time limit, demonstrating its reliability for solving real-world LP.
    \item In the case of moderate accuracy ($\epsilon=10^{-4}$), cuPDLP.jl exhibits comparable performance to Gurobi in terms of solved count and solve time, regardless of whether to use presolve. 
    For medium-sized and large-sized instances without presolve, cuPDLP.jl establishes an advantage over Gurobi, achieving a 3x speed-up on medium problems with 4 more instances solved and a 3.7x speed-up on large instances with one additional solved instance, respectively. With Gurobi presolve, cuPDLP.jl is able to solve 381 out of 383 instances and the solve time is comparable to the best of the three Gurobi methods (i.e., barrier methods).
    \item In the case of high accuracy ($\epsilon=10^{-8}$), cuPDLP.jl has comparable performance to Gurobi primal and dual simplex method, though it is inferior to the Gurobi barrier method. 
    \item Gurobi presolve can improve the performance of all Gurobi methods as well as the performance of cuPDLP.jl. The effect of presolve is more significant for Gurobi methods. This is expected because Gurobi presolve is designed to speed up Gurobi methods.
\end{itemize}
To summarize, these observations affirm that cuPDLP.jl attains comparable performance to Gurobi in \texttt{MIP Relaxations} benchmark dataset. This demonstrates that a first-order-method-based LP solver on GPU can be on par with a strong implementation of simplex and barrier methods, even in obtaining high-accuracy solutions.

\begin{figure}[ht!]
    \hspace{-1cm}
	\begin{tabular}{c c c c}
		& \includegraphics[width=0.5\textwidth]{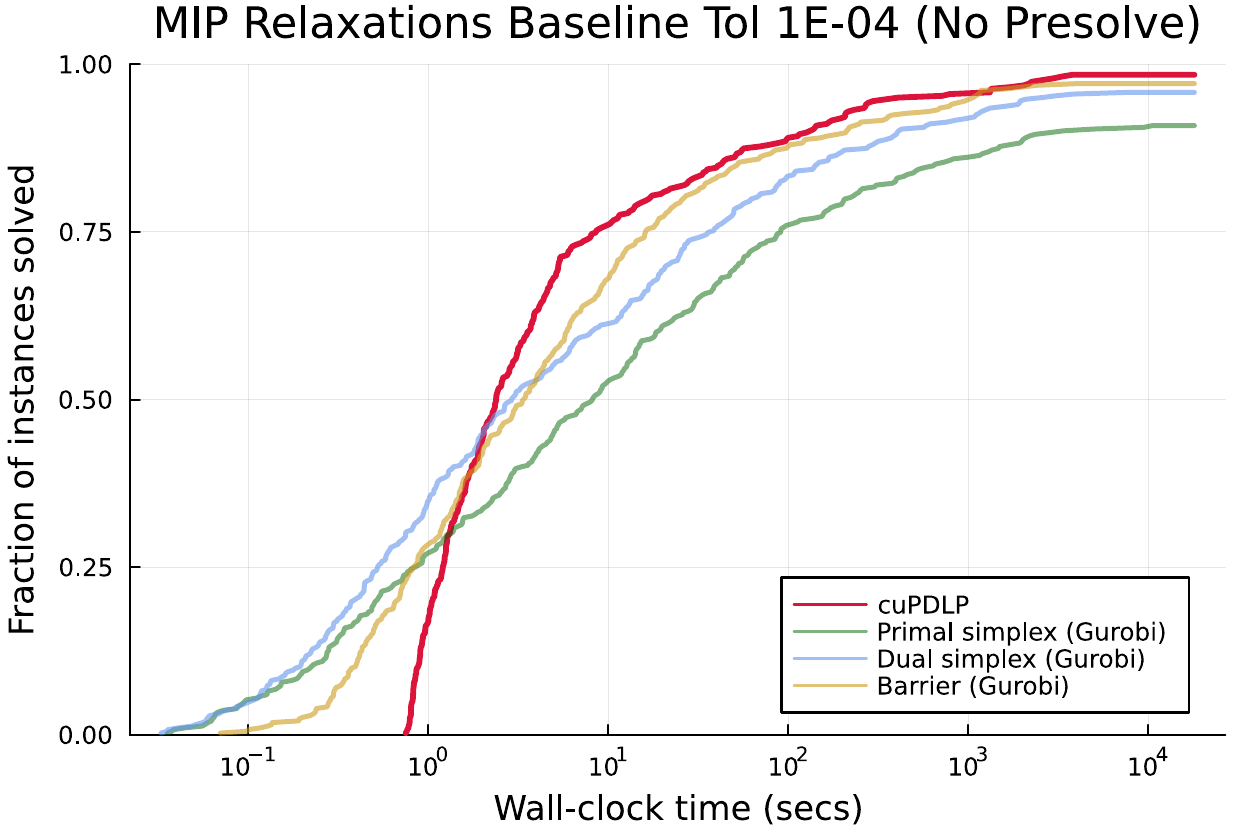}
        & \includegraphics[width=0.5\textwidth]{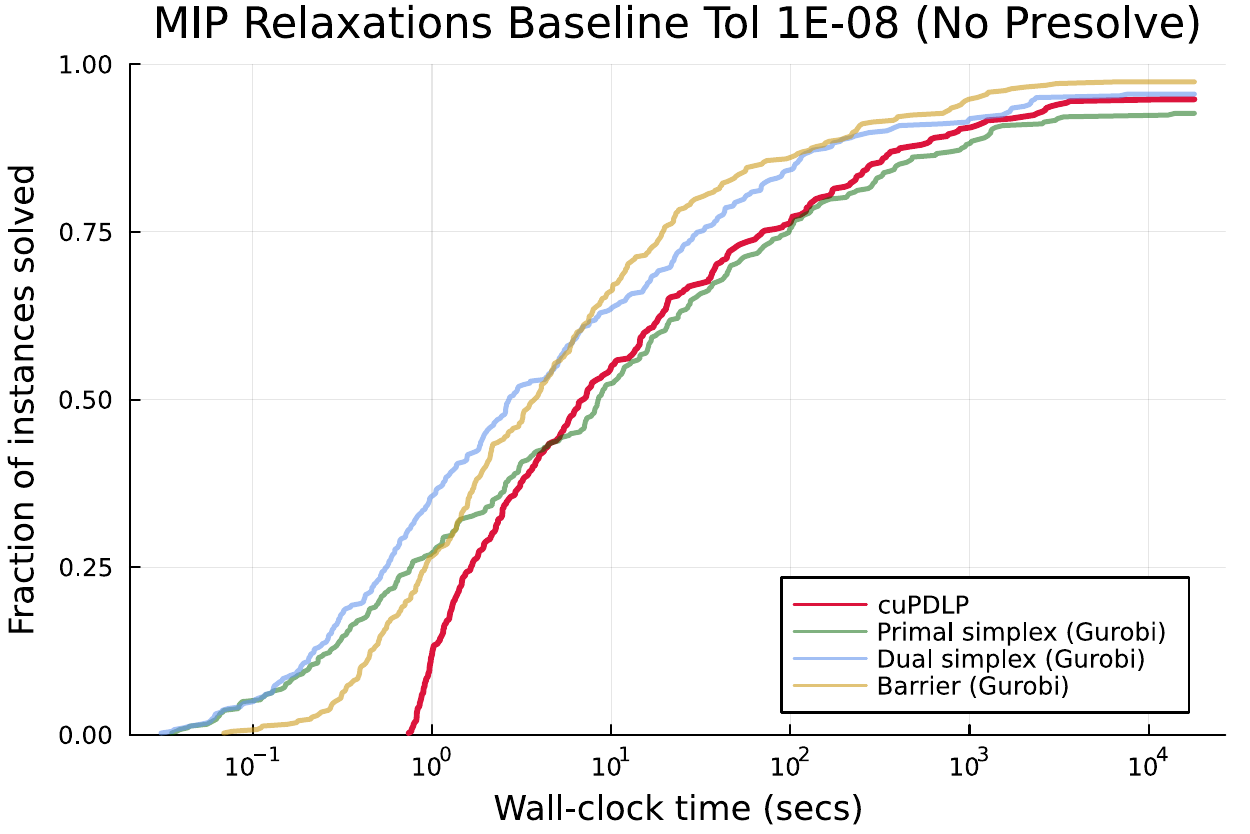}
	\end{tabular}
	\caption{Number of instances solved for \texttt{MIP Relaxations} under moderate accuracy (left) and high accuracy (right): cuPDLP.jl versus Gurobi without presolve.}
	\label{fig:performance-no-presolve}
\end{figure}

\begin{figure}[ht!]
    \hspace{-1cm}
	\begin{tabular}{c c c c}
		& \includegraphics[width=0.5\textwidth]{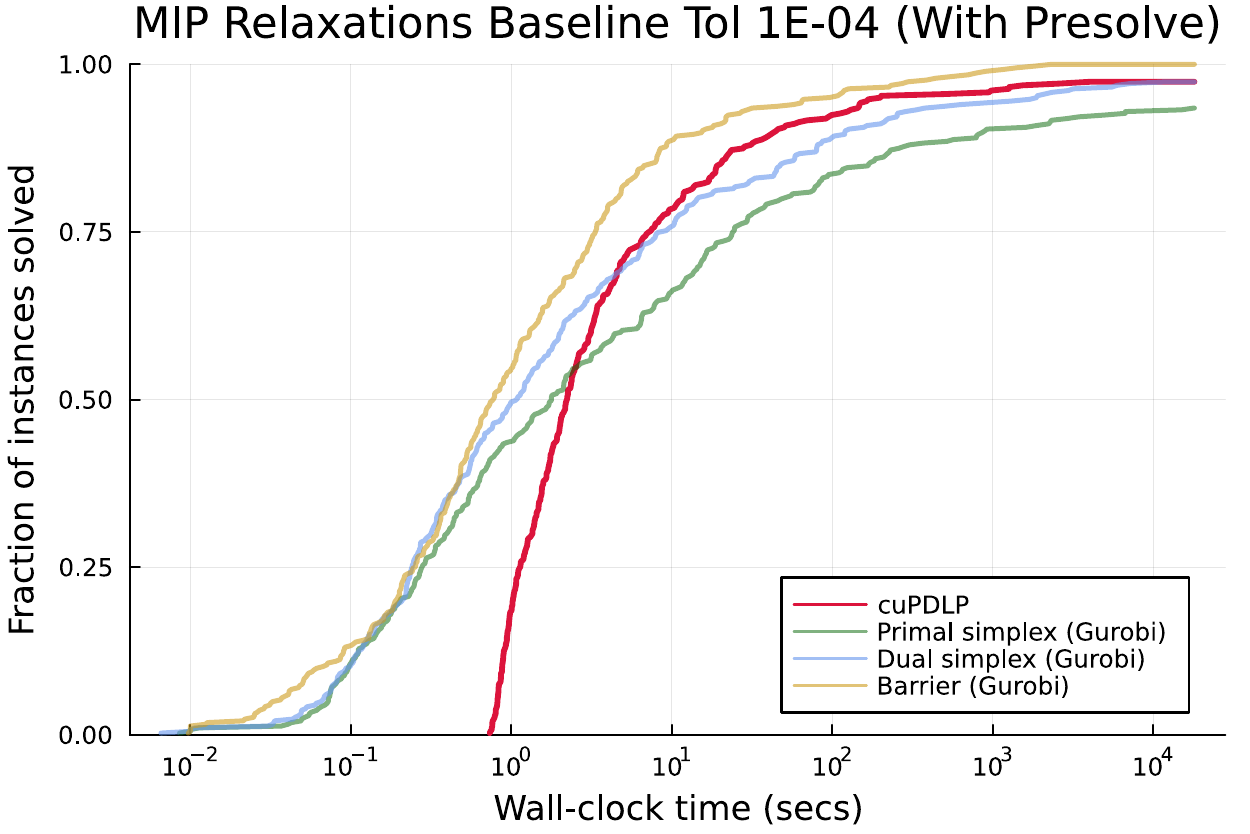}
        & \includegraphics[width=0.5\textwidth]{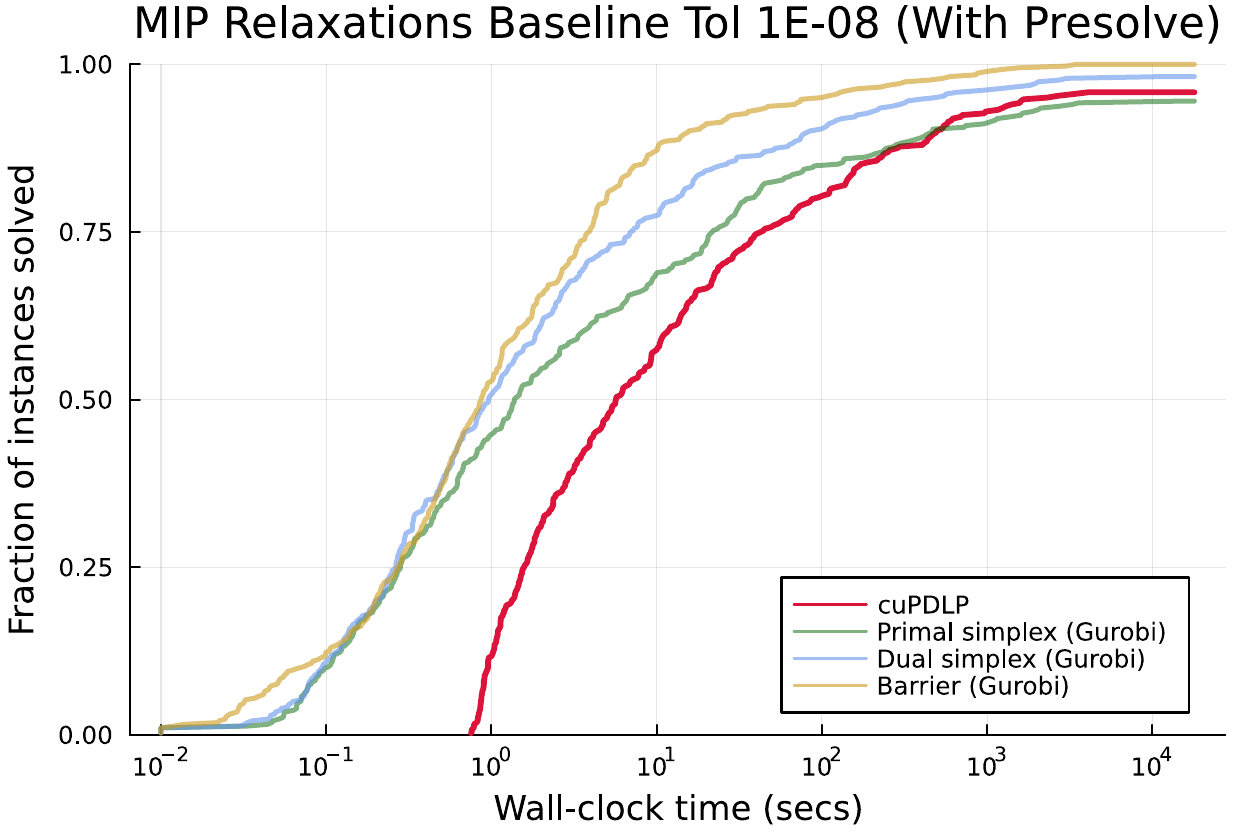}
	\end{tabular}
	\caption{Number of instances solved for \texttt{MIP Relaxations} under moderate accuracy (left) and high accuracy (right): cuPDLP.jl versus Gurobi with presolve.}
	\label{fig:performance-with-presolve}
\end{figure}

Figure \ref{fig:performance-no-presolve} and Figure \ref{fig:performance-with-presolve} show the number of solved instances of cuPDLP.jl and three methods in Gurobi on \texttt{MIP Relaxations} in a given time. The y-axes display the fraction of solved instances, and the x-axes display the wall-clock time in seconds. As shown in the left panel, when seeking solutions with moderate accuracy ($\epsilon=10^{-4}$), cuPDLP.jl has comparable performances with Gurobi barrier after 10 seconds. It eventually has better performances on \texttt{MIP Relaxation} than all three methods of Gurobi without presolve. In addition, we can see the performance of cuPDLP.jl for high-quality solution ($\epsilon=10^{-8}$), as shown in the right panel, is still comparable to Gurobi. An observation is that the number of instances Gurobi can solve for a given running time does not differ much for moderate and high accuracy; conversely, such difference is more apparent for cuPDLP.jl. This is a feature of a first-order-method-based solver. Another interesting fact is that Gurobi solves about 35\% of instances within one second, which is exactly the power of Gurobi. On the other hand, cuPDLP.jl has a computational overhead of around one second due to the GPU kernel launch time. 

\subsubsection{cuPDLP.jl versus PDLP}\label{sec:miplib-pdlp}
\begin{table}[ht!]
\centering
\begin{tabular}{cccccccccc}
\hline
\multirow{2}{*}{}                                   & \multicolumn{2}{c}{\begin{tabular}[c]{@{}c@{}}\textbf{Small (269)} \\ (1-hour limit)\end{tabular}} & \multicolumn{2}{c}{\begin{tabular}[c]{@{}c@{}}\textbf{Medium (94)}\\ (1-hour limit)\end{tabular}}    & \multicolumn{2}{c}{\begin{tabular}[c]{@{}c@{}}\textbf{\textbf{Large (20)}}\\ (5-hour limit)\end{tabular}}   & \multicolumn{2}{c}{\textbf{Total (383)}}       \\
                                                                                                   & \textbf{Count} & \textbf{Time} & \textbf{Count} & \textbf{Time} & \textbf{Count} & \textbf{Time}  & \textbf{Count} & \textbf{Time} \\ \hline
\multicolumn{1}{c}{\textbf{cuPDLP.jl}}     & 266                   & 8.61               & 92                    & 14.80               & 19                    & 111.19 &377 &12.02             \\
\multicolumn{1}{c}{\textbf{\begin{tabular}[c]{@{}c@{}}FirstOrderLp.jl\end{tabular}}}               & 253                   & 35.94               & 82                    & 155.67              & 12                    & 2002.21  & 347 & 66.67           \\
\multicolumn{1}{c}{\textbf{\begin{tabular}[c]{@{}c@{}}PDLP (1 thread)\end{tabular}}} & 256                   & 22.69               & 85                    & 98.38             & 15                    & 1622.91   &    356 & 43.81       \\
\multicolumn{1}{c}{\textbf{\begin{tabular}[c]{@{}c@{}}PDLP (4 threads)\end{tabular}}}      & 260                   & 24.03               & 91                    & 42.94               & 15                    & 736.20     & 366 & 34.57   
\\
\multicolumn{1}{c}{\textbf{\begin{tabular}[c]{@{}c@{}}PDLP (16 threads)\end{tabular}}}      & 238                   & 104.72               & 84                    & 142.79               & 15                    & 946.24     & 337 & 127.49       \\ \hline
\end{tabular}
\caption{Solve time in seconds and SGM10 of different solvers on instances of \texttt{MIP Relaxations} with tolerance $10^{-4}$: cuPDLP.jl versus PDLP.}
\label{tab:miplib-1e-4-pdlp}
\end{table}

\begin{table}[ht!]
\centering
\begin{tabular}{cccccccccc}
\hline
\multirow{2}{*}{}                                   & \multicolumn{2}{c}{\begin{tabular}[c]{@{}c@{}}\textbf{Small (269)} \\ (1-hour limit)\end{tabular}} & \multicolumn{2}{c}{\begin{tabular}[c]{@{}c@{}}\textbf{Medium (94)}\\ (1-hour limit)\end{tabular}}    & \multicolumn{2}{c}{\begin{tabular}[c]{@{}c@{}}\textbf{\textbf{Large (20)}}\\ (5-hour limit)\end{tabular}}   & \multicolumn{2}{c}{\textbf{Total (383)}}       \\
                                                                                                   & \textbf{Count} & \textbf{Time} & \textbf{Count} & \textbf{Time} & \textbf{Count} & \textbf{Time}  & \textbf{Count} & \textbf{Time} \\ \hline
\multicolumn{1}{c}{\textbf{cuPDLP.jl}}     & 261                   & 23.47                & 86                    & 40.69                & 16                    & 421.40  &363 &32.35             \\
\multicolumn{1}{c}{\textbf{FirstOrderLp.jl}}       & 235                   & 91.14                & 68                    & 389.34               & 9                     & 3552.50    &312 &160.63          \\
\multicolumn{1}{c}{\textbf{\begin{tabular}[c]{@{}c@{}}PDLP (1 thread)\end{tabular}}} & 250                   & 49.31               & 73                    & 259.04             & 12                    & 3818.42   &    335 & 96.86       \\
\multicolumn{1}{c}{\textbf{\begin{tabular}[c]{@{}c@{}}PDLP (4 threads)\end{tabular}}}      & 245                   & 54.19               & 81                   & 136.16               & 14                    & 1789.54     & 340 & 83.49         \\ 
\multicolumn{1}{c}{\textbf{\begin{tabular}[c]{@{}c@{}}PDLP (16 threads)\end{tabular}}}      & 214                   & 248.34               & 69                   & 403.17               & 14                    & 2475.57     & 297 & 316.27        \\ 
\hline
\end{tabular}
\caption{Solve time in seconds and SGM10 of different solvers on instances of \texttt{MIP Relaxations} with tolerance $10^{-8}$: cuPDLP.jl versus PDLP.}
\label{tab:miplib-1e-8-pdlp}
\end{table}
Tables \ref{tab:miplib-1e-4-pdlp} and \ref{tab:miplib-1e-8-pdlp} present the performance comparison of cuPDLP.jl and its CPU implementations on \texttt{MIP Relaxations} with tolerances of $10^{-4}$ and $10^{-8}$, respectively. The two tables demonstrate that GPU can significantly speed up PDLP, in particular for large instances:
\begin{itemize}
    \item 
    For moderate accuracy (Table \ref{tab:miplib-1e-4-pdlp}), cuPDLP.jl demonstrates a 4x speed-up for small instances, a 10x speed-up for medium instances, and a 20x speed-up for large instances, compared to FirstOrderLp.jl, a CPU implementation of PDLP in Julia. When comparing cuPDLP.jl with the more delicate C++ implementation PDLP with multithreading support, it still exhibits a significant speedup for all sizes of problems. The significant speed-up can also be observed for high accuracy (Table \ref{tab:miplib-1e-8-pdlp}).
    \item In terms of solved count, cuPDLP.jl solves significantly more instances regardless of the scales. In particular, comparing with FirstOrderLp.jl under tolerance $\epsilon=10^{-4}$, cuPDLP.jl solves 13 more small-sized instances, 10 more medium-sized instances and 7 more large problems, with in total 30 more instances solved. Compared to PDLP with the best of 1 thread, 4 threads or 16 threads, cuPDLP.jl can solve 6 more small-sized instances and 4 more large-sized instances. The improvement is also remarkable when looking at results for high accuracy $\epsilon=10^{-8}$.
\end{itemize}
In summary, the GPU-implemented cuPDLP.jl consistently outperforms the CPU-implemented PDLP in numerical performance on \texttt{MIP Relaxations}, with cuPDLP.jl demonstrating even more pronounced advantages on instances of medium to large scale.

\begin{figure}[ht!]
    \hspace{-1cm}
	\begin{tabular}{c c c c}
		& \includegraphics[width=0.5\textwidth]{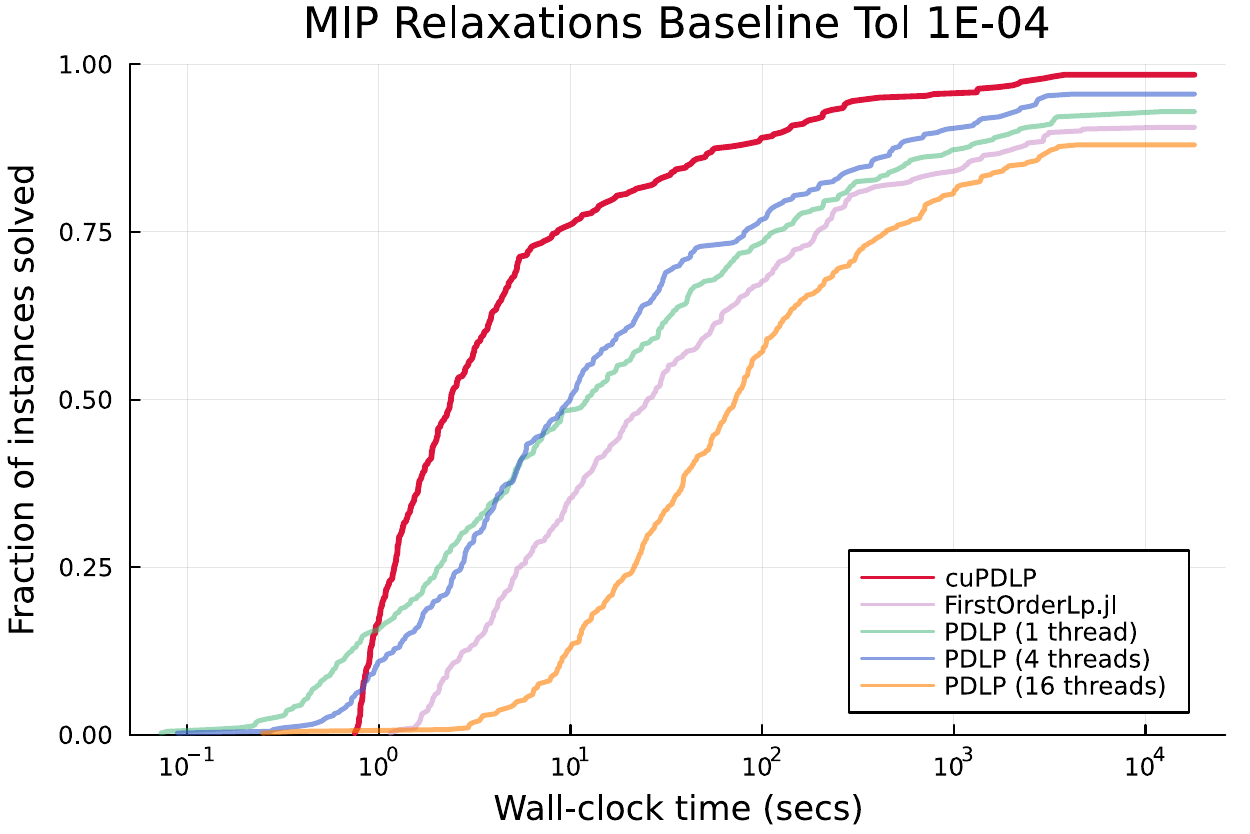}
        & \includegraphics[width=0.5\textwidth]{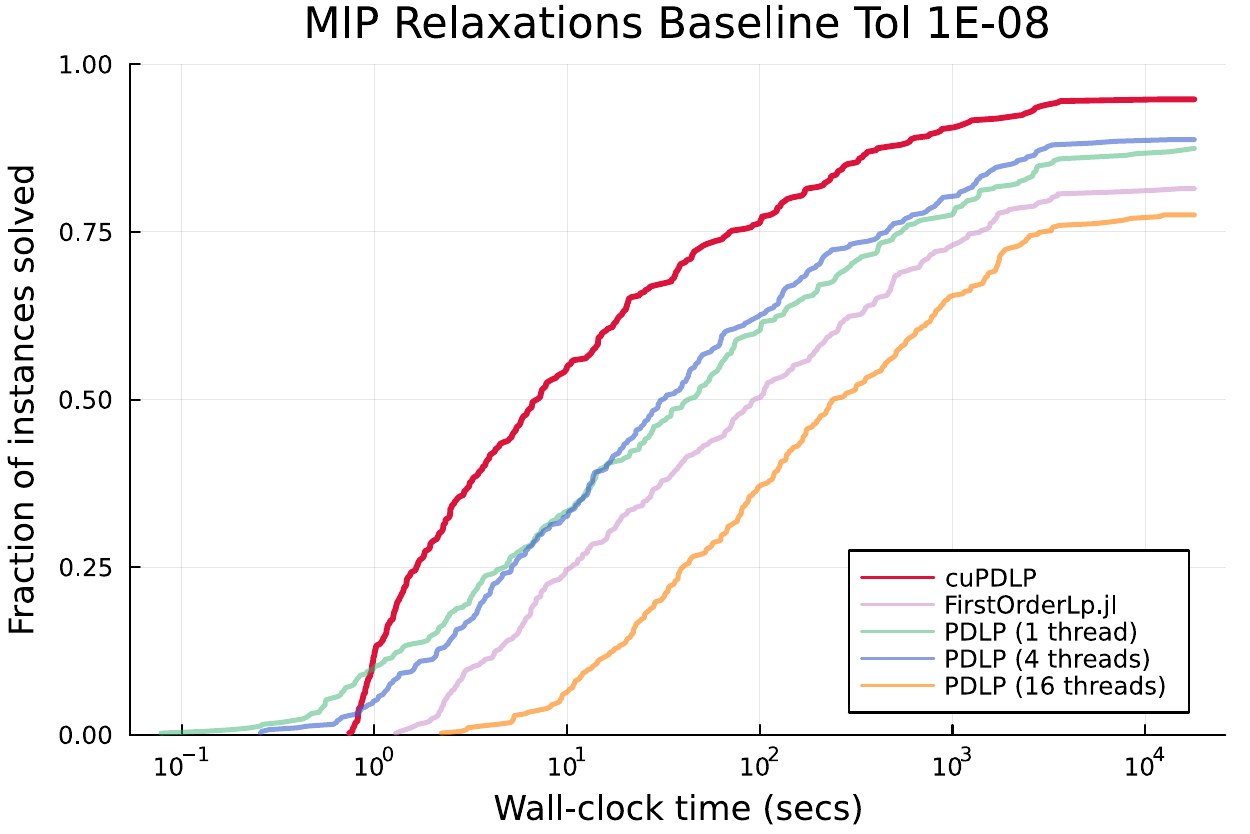}
	\end{tabular}
	\caption{Number of instances solved for \texttt{MIP Relaxations} under moderate accuracy (left) and high accuracy (right): cuPDLP.jl versus PDLP.}
	\label{fig:performance-pdlp}
\end{figure}

Similar to Figure \ref{fig:performance-no-presolve} and \ref{fig:performance-with-presolve}, Figure \ref{fig:performance-pdlp} demonstrates the number of solved instances of cuPDLP.jl and three CPU implementations of PDLP on \texttt{MIP Relaxations} in a given time. As shown in both panels, when seeking solutions with moderate accuracy ($\epsilon=10^{-4}$) and high accuracy ($\epsilon=10^{-8}$), cuPDLP.jl has a clear superior performance to all CPU versions of PDLP. Moreover, it is notable that cuPDLP.jl has a computational overhead of around one second due to the GPU kernel launch time.

\subsection{Mittelmann's LP benchmark set}\label{sec:mittelmann}
We further compare cuPDLP.jl with PDLP and Gurobi on \texttt{Mittelmann's LP} benchmark set. \texttt{Mittelmann's LP} benchmark set~\cite{mittelmann} is a standard LP benchmark to test the numerical performances of different LP solvers and includes 49 LP instances in our experiments. 

Table \ref{tab:mittelmann-no-presolve}, Table \ref{tab:mittelmann-with-presolve} and Table \ref{tab:mittelmann-pdlp} summarize the results on \texttt{Mittelmann's LP} benchmark set. The observations are consistent with the findings on \texttt{MIP Relaxations} as discussed in Section \ref{sec:miplib}. Consider the moderate accuracy. cuPDLP.jl can solve more instances against all versions of PDLP and have comparable performance to Gurobi. When solving for high-accuracy solutions, the performance of cuPDLP.jl still performs inferior to Gurobi barrier but comparable to primal and dual simplex. This again demonstrates the robust performance of GPU-implemented cuPDLP.jl.

\begin{table}[ht!]
\centering
\begin{tabular}{ccccc}
\hline
\multicolumn{1}{l}{}                                & \multicolumn{2}{c}{\textbf{Tol 1E-04}}      & \multicolumn{2}{c}{\textbf{Tol 1E-08}}       \\ 
                                                    & \textbf{Count} & \textbf{Time} & \textbf{Count} & \textbf{Time} \\ \hline
\multicolumn{1}{c}{\textbf{cuPDLP.jl}}     & 44                     &   71.26             &  40                   &    231.91           \\
\multicolumn{1}{c}{\textbf{Primal simplex (Gurobi)}} & 35                    & 1937.78              & 34                   & 1715.62              \\
\multicolumn{1}{c}{\textbf{Dual simplex (Gurobi)}}   & 36                    & 1201.48              & 37                    & 1116.68              \\
\multicolumn{1}{c}{\textbf{Barrier (Gurobi)}}        & 45                    & 108.29              & 44                    & 127.58              \\ \hline
\end{tabular}
\caption{Solve time in seconds and SGM10 of different solvers on instances of \texttt{Mittelmann's LP} benchmark set with tolerance $10^{-4}$ and $10^{-8}$ without presolve.}
\label{tab:mittelmann-no-presolve}
\end{table}

\begin{table}[ht!]
\centering
\begin{tabular}{ccccc}
\hline
\multicolumn{1}{l}{}                                & \multicolumn{2}{c}{\textbf{Tol 1E-04}}      & \multicolumn{2}{c}{\textbf{Tol 1E-08}}       \\ 
                                                    & \textbf{Count} & \textbf{Time} & \textbf{Count} & \textbf{Time} \\ \hline
\multicolumn{1}{c}{\textbf{cuPDLP.jl}}     & 47                    &  52.44              & 44                    & 167.15             \\
\multicolumn{1}{c}{\textbf{Primal simplex (Gurobi)}} & 42                    & 606.03              & 39                    & 499.09              \\
\multicolumn{1}{c}{\textbf{Dual simplex (Gurobi)}}   & 40                    & 290.69              & 42                    & 212.00              \\
\multicolumn{1}{c}{\textbf{Barrier (Gurobi)}}        & 48                    & 24.16              & 47                    & 34.18              \\ \hline
\end{tabular}
\caption{Solve time in seconds and SGM10 of different solvers on instances of \texttt{Mittelmann's LP} benchmark set with tolerance $10^{-4}$ and $10^{-8}$ with presolve.}
\label{tab:mittelmann-with-presolve}
\end{table}

\begin{table}[ht!]
\centering
\begin{tabular}{ccccc}
\hline
\multicolumn{1}{l}{}                                & \multicolumn{2}{c}{\textbf{Tol 1E-04}}      & \multicolumn{2}{c}{\textbf{Tol 1E-08}}       \\ 
                                                    & \textbf{Count} & \textbf{Time} & \textbf{Count} & \textbf{Time} \\ \hline
\multicolumn{1}{c}{\textbf{cuPDLP.jl}}     & 44                     &   71.26             &  40                   &    231.91           \\
\multicolumn{1}{c}{\textbf{FirstOrderLp.jl}}       & 34                    & 917.49             & 25                    & 2504.79             \\
\multicolumn{1}{c}{\textbf{PDLP (1 thread)}} & 39                    & 586.50             & 31                    & 1661.20            \\
\multicolumn{1}{c}{\textbf{PDLP (4 threads)}}   & 40                    & 302.54              & 34                    & 930.41                        \\ 
\multicolumn{1}{c}{\textbf{PDLP (16 threads)}}   & 39                    & 776.22              & 29                    & 2171.95                        \\ \hline
\end{tabular}
\caption{Solve time in seconds and SGM10 of different solvers on instances of \texttt{Mittelmann's LP} benchmark set with tolerance $10^{-4}$ and $10^{-8}$.}
\label{tab:mittelmann-pdlp}
\end{table}

\section{Conclusion}
In this paper, we present cuPDLP.jl, a GPU implementation of restarted PDHG for solving LP in Julia. The numerical experiments demonstrate that the prototype GPU implementation cuPDLP.jl can have comparable performance with commercial solvers like Gurobi and superior performance on large instances. This sheds light on using GPU to develop high-performance optimization solvers.

\section*{Acknowledgement}
The authors would like to thank Azam Asl and Miles Lubin for the early discussions on the GPU implementation of PDLP, and David Applegate for his encouragement and support.

\bibliographystyle{amsplain}
\bibliography{ref-papers}

\appendix
\section{Theoretical guarantees of restarted PDHG with KKT error for LP}\label{app:thoery}
Consider LP of the standard form:
\begin{equation}\label{eq:lp-standard}
    \begin{aligned}[c]
    \min_{x\in \mathbb R^n}~ c^\top x \ ~~
\text{s.t.}\ Ax = b,\ x \geq 0\ .
    \end{aligned}
\end{equation} 
The KKT error at solution $z=(x,y)$ is defined as the norm of the violation of KKT system of \eqref{eq:lp-standard}
{\small
\begin{equation*}
    \mathrm{KKT}(z)=\mathrm{KKT}(x,y)=\left\|  \begin{pmatrix}
        Ax-b \\ [-x]^+ \\ [A^\top y-c]^+ \\ [c^\top x-b^\top y]^+
    \end{pmatrix}\right\|_2 \ ,
\end{equation*}
}
where $[x]^+= [(x_1,...,x_n)]^+=(\max\{x_1,0\},\ldots,\max\{x_n,0\})$ is the positive part of the vector $x$.

\begin{algorithm}[H]
\caption{Restarted PDHG for solving \eqref{eq:lp-standard}}
\label{alg:pdhg-restart-kkt}
\SetKwInOut{Input}{Input}
\Input{Initial point $(x^0,y^0)$, outer loop counter $n\leftarrow 0$. Step-sizes $\tau$, $\sigma$. Restart decay $\beta\in (0,1)$.}

\Repeat{\upshape $(x^{n+1,0},y^{n+1,0})$ convergence}{
  initialize the inner loop counter $k\leftarrow0$;\\
  \Repeat{$\mathrm{KKT}(\bar x^{n,k+1},\bar y^{n,k+1})\leq \beta \mathrm{KKT}(x^{n,0}, y^{n,0})$ \upshape holds}{
    $x^{n,k+1}\leftarrow \text{proj}_{\mathbb R_+^n}(x^{n,k}-\tau (c-A^\top y^{n,k}))$;\\
    $y^{n,k+1}\leftarrow y^{n,k}+\sigma (b-A(2x^{n,k+1}-x^{n,k}))$;\\
    $\bar x^{n,k+1}\leftarrow \frac{k}{k+1}\bar x^{n,k}+\frac{1}{k+1}x^{n,k+1}$;\\
    $\bar y^{n,k+1}\leftarrow \frac{k}{k+1}\bar y^{n,k}+\frac{1}{k+1}y^{n,k+1}$;
  }
  initialize the initial solution $(x^{n+1,0}, y^{n+1,0})\leftarrow(\bar x^{n,k+1}, \bar y^{n,k+1})$;\\
  $n\leftarrow n+1$;
}
\end{algorithm}

Consider restarted PDHG with KKT error as restarting criteria (Algorithm \ref{alg:pdhg-restart-kkt}).
Denote $z^{n,k}=(x^{n,k},y^{n,k})$ (and $\bar z^{n,k}=(\bar x^{n,k},\bar y^{n,k})$, correspondingly) as the primal-dual solution pair at the $n$-th other iteration and the $k$-th inner iteration. Without loss of generality, assume the primal and dual stepsizes are equal, i.e., $\tau=\sigma=:s$ (otherwise we can rescale the primal and the dual problem so that they share the same step-size; a similar proof strategy is used in \cite{applegate2023faster}). Denote $P_s=\begin{pmatrix}
    I & sA^\top \\ sA & I
\end{pmatrix}$ and define $\|v\|_{P_s}^2=\langle v, P_sv\rangle$. Denote $\|\cdot\|_2$ the Euclidean norm and $\mathrm{dist}_2(w,U)=\min_{u\in U}\|w-u\|_2$ the distance between point $w$ and set $U$.

First, it is straight-forward to see that $\mathrm{KKT}$ is a sharp function, i.e.,

\begin{prop}[{\cite{hoffman1952approximate,pena2021new}}]
    There exists a constant $\alpha>0$ such that for any $z=(x,y)$, it holds that
    \begin{equation*}
        \alpha\mathrm{dist}_2(z,\mathcal Z^*)\leq \mathrm{KKT}(z)\ .
    \end{equation*}
\end{prop}

The next theorem proves the linear convergence rate of Algorithm \ref{alg:pdhg-restart-kkt}. Notice that this is essentially the same linear convergence rate as that in \cite{applegate2023faster} up to a constant.
\begin{thm}
Consider $\{z^{n,k}\}$ the iterates of Algorithm \ref{alg:pdhg-restart-kkt} for solving \eqref{eq:primal-dual}. Suppose stepsize $s\leq \frac{1}{2\|A\|_2}$. Then it holds that
    \begin{enumerate}
    \item[(i)] There exists an $R>0$ such that for any $n$ and $k$, $\|z^{n,k}\|_2\leq R$.
    \item[(ii)] The restart length $\tau^n$ is upper bounded by
    \begin{equation*}
        \tau^n\leq \frac{8\sqrt 2\sqrt{1+R^2}}{\alpha s\beta} \ .
    \end{equation*}
    \item[(iii)] The distance to optimal solution set decays linearly
    \begin{equation*}
        \mathrm{dist}(z^{n,0},\mathcal Z^*) \leq \beta^{n+1}\frac{1}{\alpha}\mathrm{KKT}(z^{0,0}) \ .
    \end{equation*}
\end{enumerate}
\end{thm}
\begin{rem}
    Let stepsize $s=\frac{C}{\|A\|_2}$ where constant $C\leq \frac{1}{2}$. The iteration complexity of Algorithm \ref{alg:pdhg-restart-kkt} to find an $\epsilon$-solution equals $\widetilde O\pran{\frac{\|A\|_2}{\alpha}\log\frac{1}{\epsilon}}$, which is optimal (upto log terms) in the sense of matching the lower complexity bounds~\cite[Corollary 1]{applegate2023faster}.
\end{rem}

\begin{proof}
(i) Let $R=\sqrt 2 \|z^{0,0}-z^*\|_{P_s} + \|z^*\|_2$. By non-expansiveness of PDHG~\cite{chambolle2016ergodic}, we have
{\footnotesize
\begin{equation*}
    \|z^{n,k}-z^*\|_{P_s}\leq \|z^{n,0}-z^*\|_{P_s}=\left\|\tfrac{1}{\tau^{n-1}}\sum_{k=1}^{\tau^{n-1}}z^{n-1,k}-z^*\right\|_{P_s}\leq \frac{1}{\tau^{n-1}}\sum_{k=1}^{\tau^{n-1}}\left\|z^{n-1,k}-z^*\right\|_{P_s}\leq \|z^{n-1,0}-z^*\|_{P_s} \leq \ldots\leq  \|z^{0,0}-z^*\|_{P_s} \ .
\end{equation*}
}
Thus we achieve
\begin{equation*}
    \|z^{n,k}\|_2 \leq \|z^{n,k}-z^*\|_{2} +\|z^*\| \leq \sqrt 2 \|z^{0,0}-z^*\|_{P_s} + \|z^*\|_2=R \ .
\end{equation*}

(ii) Suppose $\tau^n> \frac{8\sqrt 2\sqrt{1+R^2}}{\alpha\beta}=:k^*$. Note that
\begin{equation*}
\begin{aligned}
    \mathrm{KKT}(\bar z^{n-1,k^*})& \leq \sqrt{1+R^2}\rho_{\|\bar z^{n-1,k^*}-z^{n-1,0}\|_{P_s}}(\bar z^{n-1,k^*})\leq \frac{4\sqrt{1+R^2}}{s k^*}\|\bar z^{n-1,k^*}-z^{n-1,0}\|_{P_s}\\
    & \leq \frac{8\sqrt{1+R^2}}{s k^*}\mathrm{dist}_{P_s}(z^{n-1,0},\mathcal Z^*)\leq \frac{8\sqrt 2\sqrt{1+R^2}}{s k^*}\mathrm{dist}_{2}(z^{n-1,0},\mathcal Z^*) \leq \frac{8\sqrt 2\sqrt{1+R^2}}{\alpha s k^*}\mathrm{KKT}(z^{n-1,0})\\
    & \leq \beta \mathrm{KKT}(z^{n-1,0}) \ ,
\end{aligned}
\end{equation*}
where the first and fifth inequalities follow from \cite[Lemma 4]{applegate2023faster}, while the second and third inequalities utilize \cite[Property 3]{applegate2023faster}. The fourth inequality uses $P_s\preceq 2I$ which is a direct consequence of stepsize $s\leq \frac{1}{2\|A\|_2}$. The last inequality follows from the definition of $k^*=\frac{8\sqrt 2\sqrt{1+R^2}}{\alpha\beta}$.

This implies the restart condition holds at $k^*$. We finish the proof by noticing that this contradicts $\tau^n>k^*$.

(iii) Note that the KKT error has linearly decay due to the adaptive restart scheme
    \begin{equation*}
        \begin{aligned}
            \mathrm{KKT}(z^{n,0})=\mathrm{KKT}(\bar z^{n-1,\tau^{n-1}})\leq \beta \mathrm{KKT}(z^{n-1,0})\leq \ldots \leq \beta^{n} \mathrm{KKT}(z^{0,0}) \ ,
        \end{aligned}
    \end{equation*}
    and thus it follows from sharpness of KKT error that
    \begin{equation*}
        \mathrm{dist}_2(z^{n,0},\mathcal Z^*)\leq \frac{1}{\alpha}\mathrm{KKT}(z^{n,0})\leq \frac{1}{\alpha}\beta^{n} \mathrm{KKT}(z^{0,0}) \ .
    \end{equation*}
\end{proof}

\section{Comparison with GPU-based SCS}\label{app:comparison}
\begin{table}[ht!]
\centering
\begin{tabular}{lll}
\hline
                      & \textbf{Count} & \textbf{Time}    \\ \hline
\textbf{SCS-GPU (Tol=1E-04)}   & 28    & 1905.00 \\
\textbf{cuPDLP.jl (Tol=1E-04)} & 92    & 14.80   \\
\textbf{cuPDLP.jl (Tol=1E-08)} & 86    & 40.69   \\ \hline
\end{tabular}
\caption{Solve time in seconds and SGM10 of different solvers on medium-sized instances of \texttt{MIP Relaxations}: cuPDLP.jl versus SCS with GPU linear system solver.}
\label{tab:scs}
\end{table}
In this section, we compare the results of cuPDLP.jl with GPU-based SCS which solves the linear system in each iteration via conjugate gradient method implemented on GPU. We select medium-sized instances (94 in total) in \texttt{MIP Relaxations} as the test set. Table \ref{tab:scs} presents the results of GPU-based SCS under moderate tolerance $10^{-4}$ and cuPDLP.jl under both $10^{-4}$ and $10^{-8}$ tolerances. The advantage of cuPDLP over SCS is quite significant: GPU-based SCS only solves 28 out of 94 instances under moderate accuracy. In comparison, cuPDLP solves 92 out of 94 instances under same tolerance and 86 out of 94 instances under even higher tolerance. This showcases cuPDLP, the LP-specialized FOM-based solver, performs better at solving LP instances than general FOM-based conic solver SCS on GPU. This result is not surprising because SCS targets a more general class of optimization problems, and the CPU implementation of PDLP has showcased superior performance than CPU-based SCS~\cite{applegate2021practical}.  

\section{Figures}\label{app:figure}
\begin{figure}[ht!]
    \hspace{-1cm}
	\begin{tabular}{c c c c}
		& \includegraphics[width=0.33\textwidth]{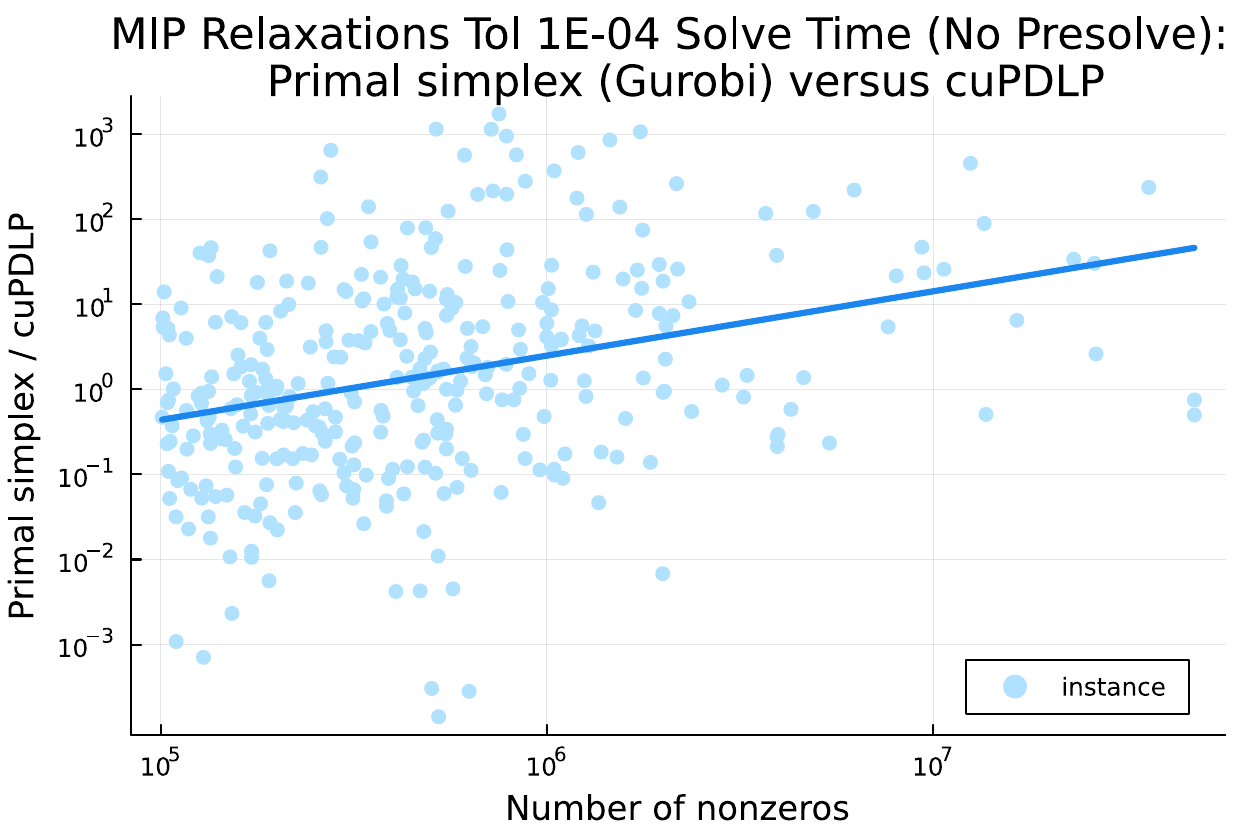}
        & \includegraphics[width=0.33\textwidth]{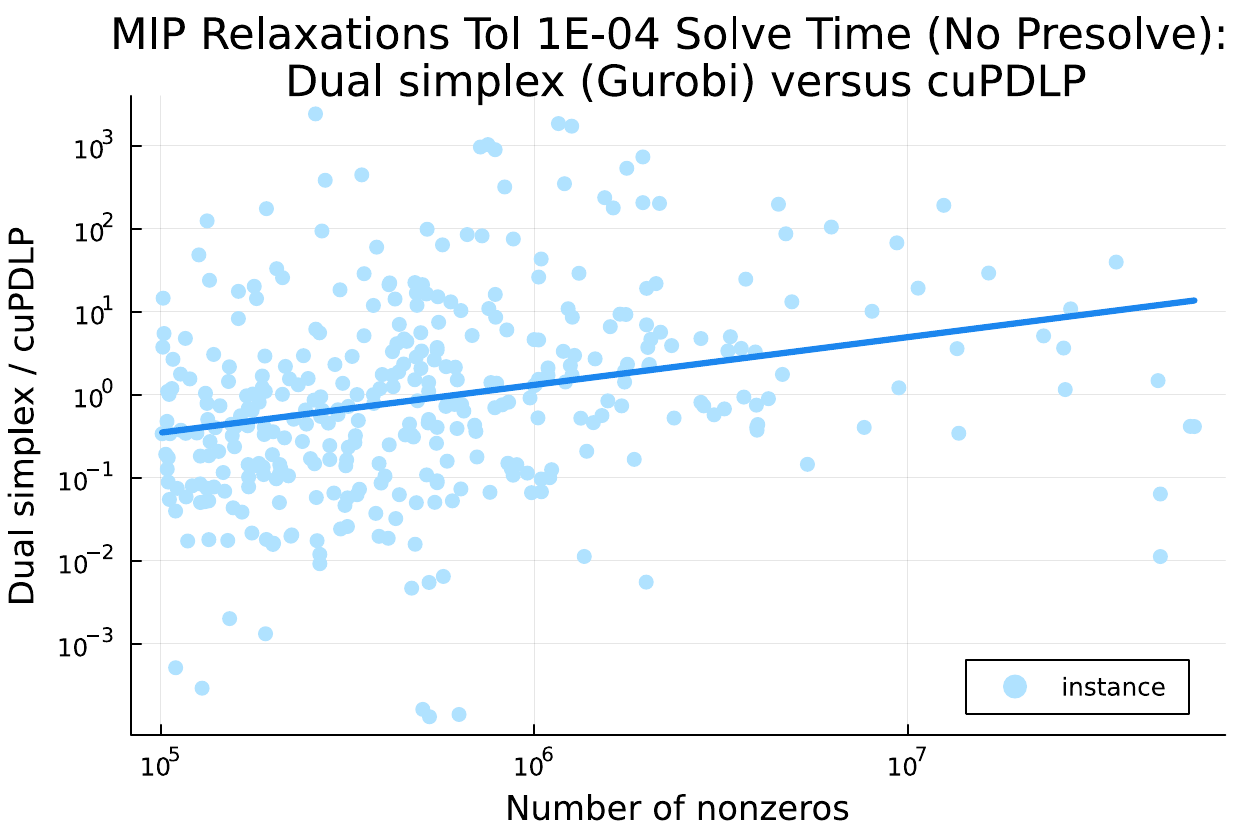}
        & \includegraphics[width=0.33\textwidth]{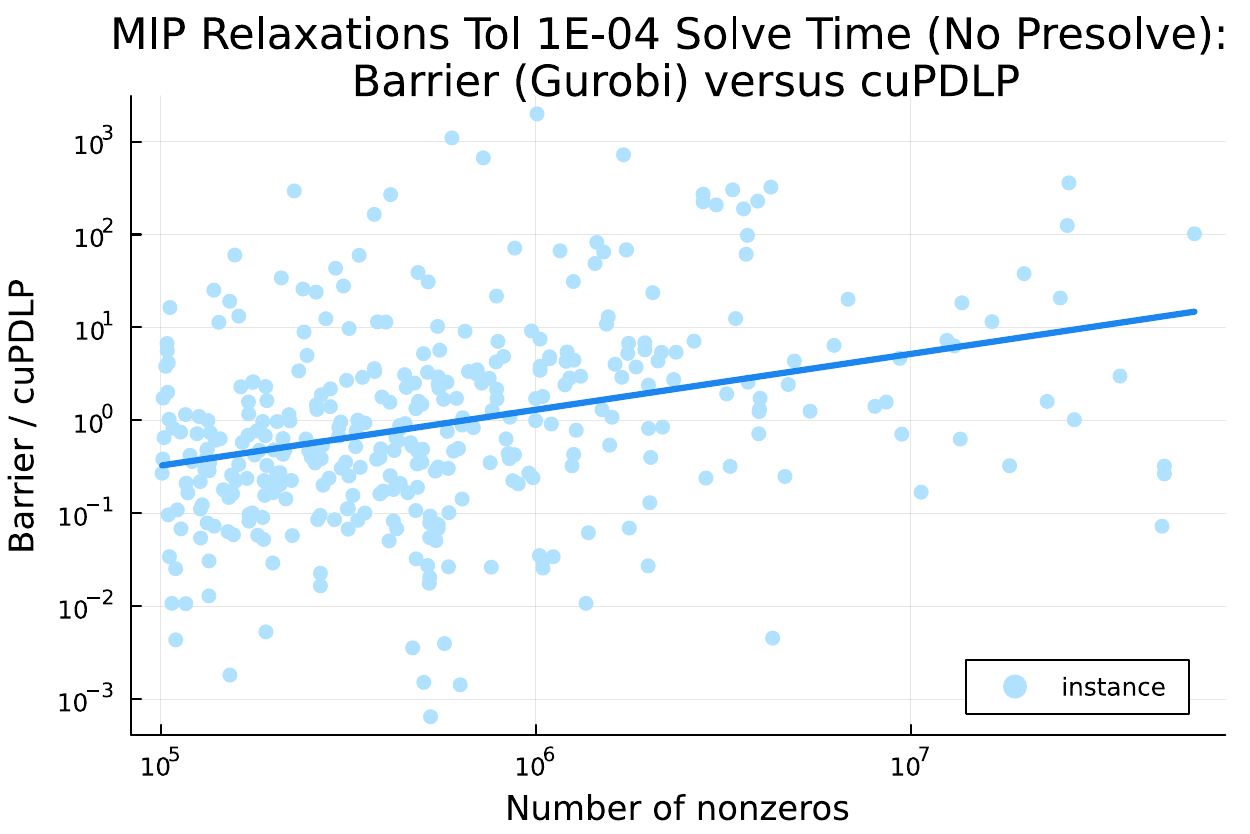}\\
        & \includegraphics[width=0.33\textwidth]{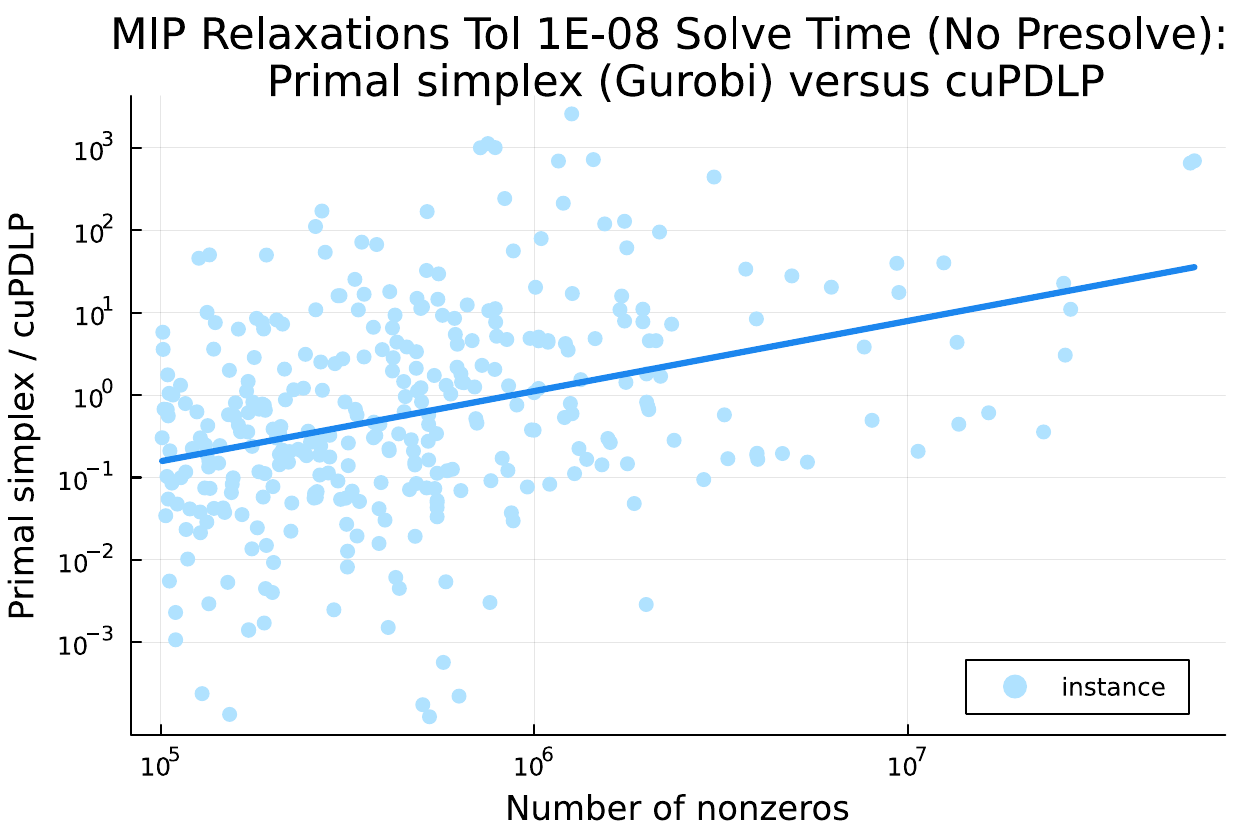}
        & \includegraphics[width=0.33\textwidth]{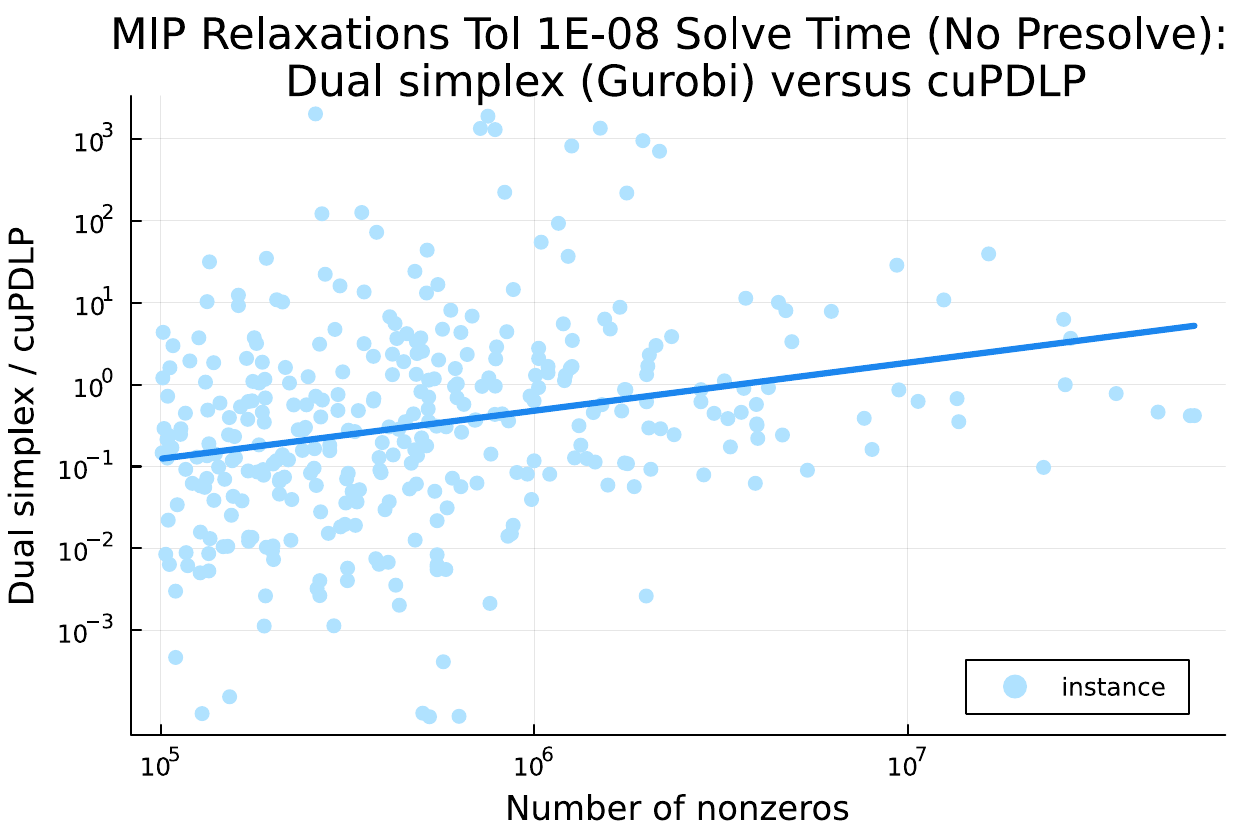}
        & \includegraphics[width=0.33\textwidth]{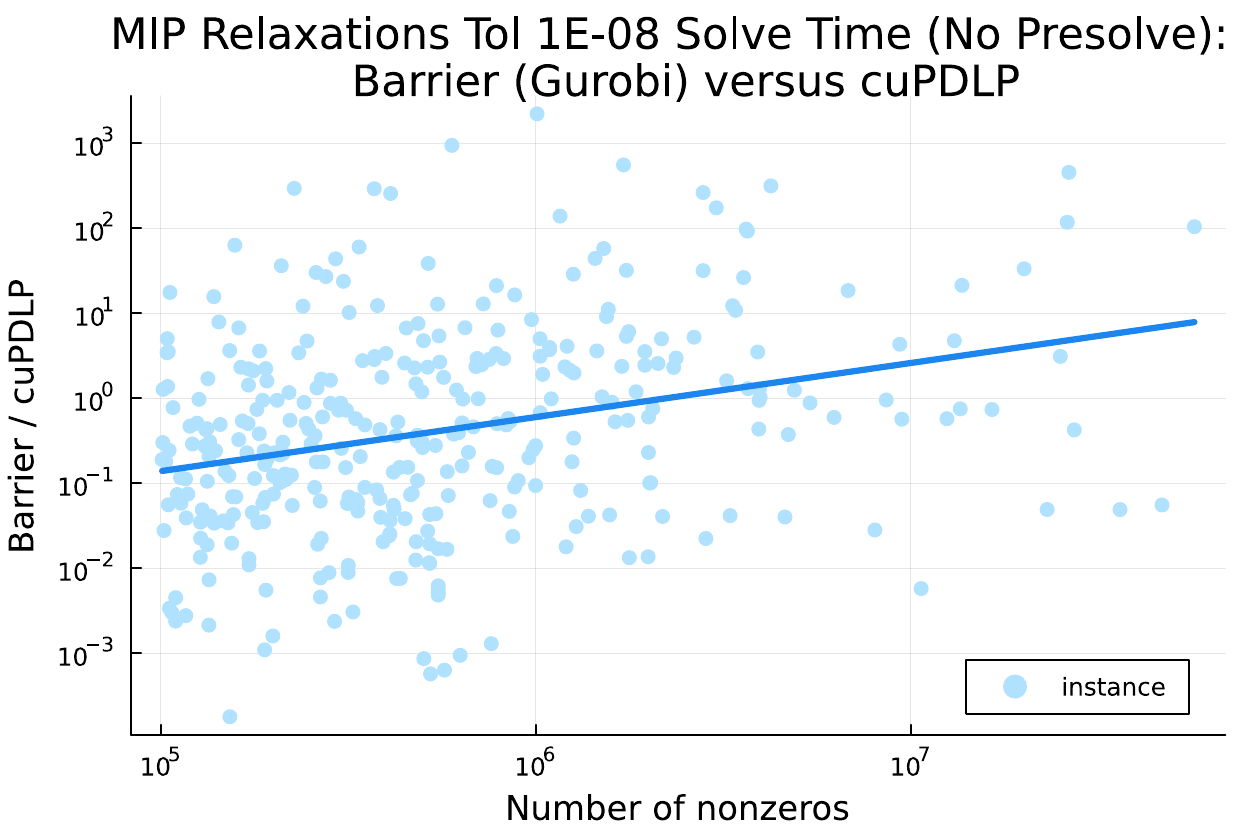}
	\end{tabular}
	\caption{Ratio of Gurobi solve time over cuPDLP.jl solve time for moderate accuracy (top) and high accuracy (bottom). Presolve is not used.}
	\label{fig:scatterplot-gurobi-no-presolve}
\end{figure}

\begin{figure}[ht!]
    \hspace{-1cm}
	\begin{tabular}{c c c c}
		& \includegraphics[width=0.33\textwidth]{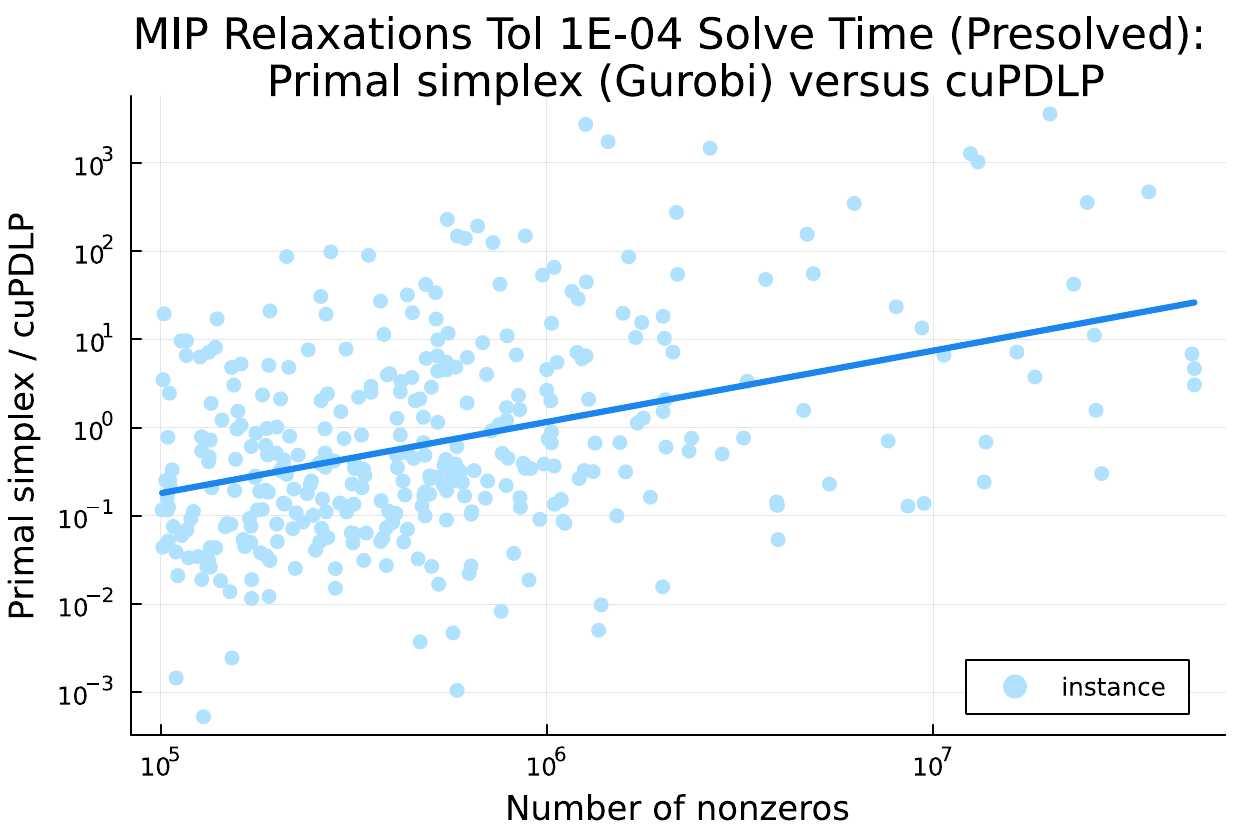}
        & \includegraphics[width=0.33\textwidth]{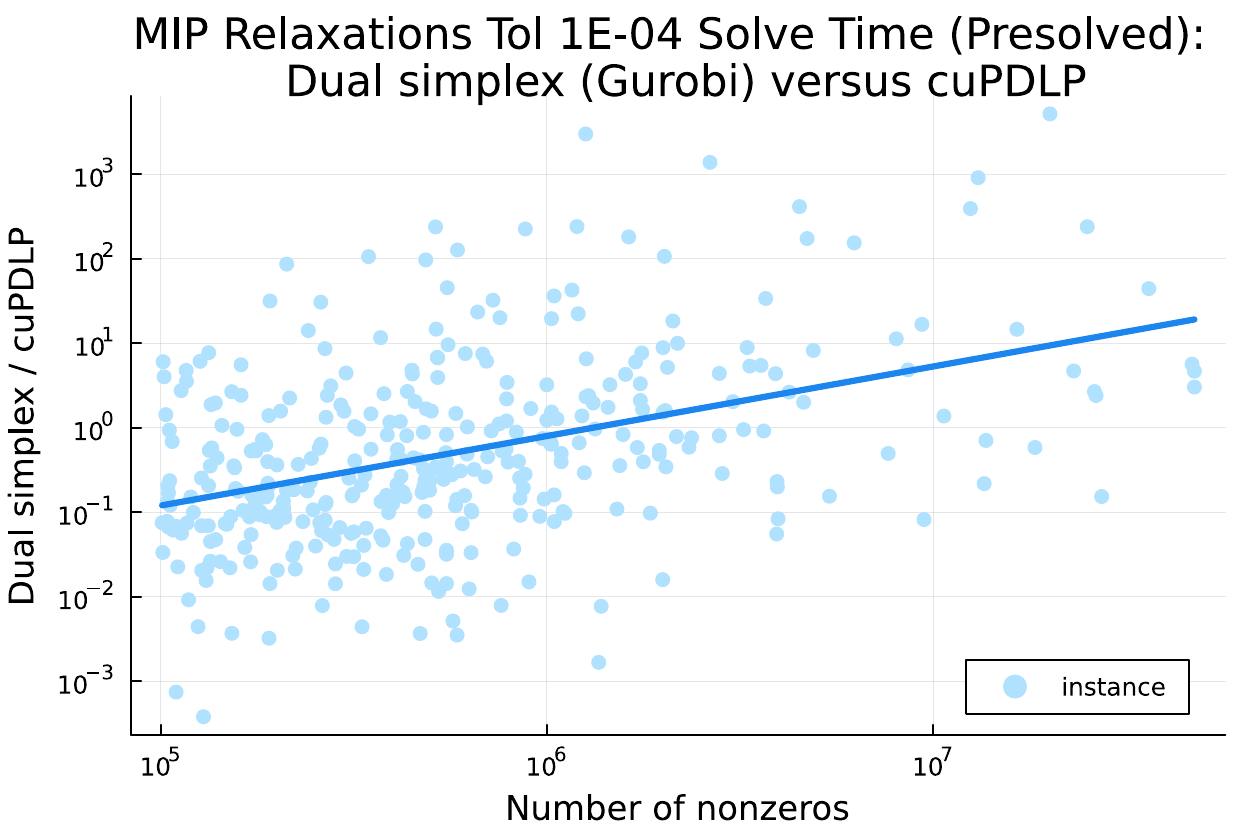}
        & \includegraphics[width=0.33\textwidth]{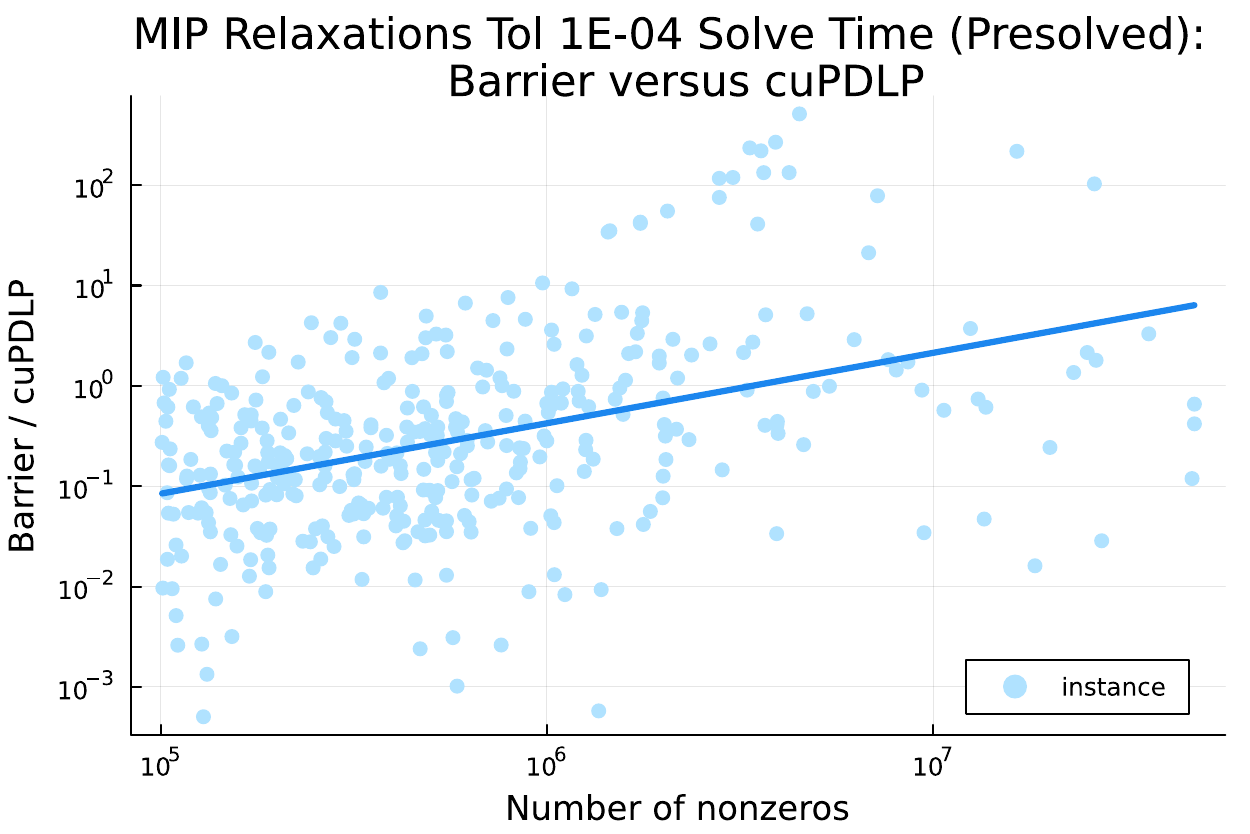}\\
        & \includegraphics[width=0.33\textwidth]{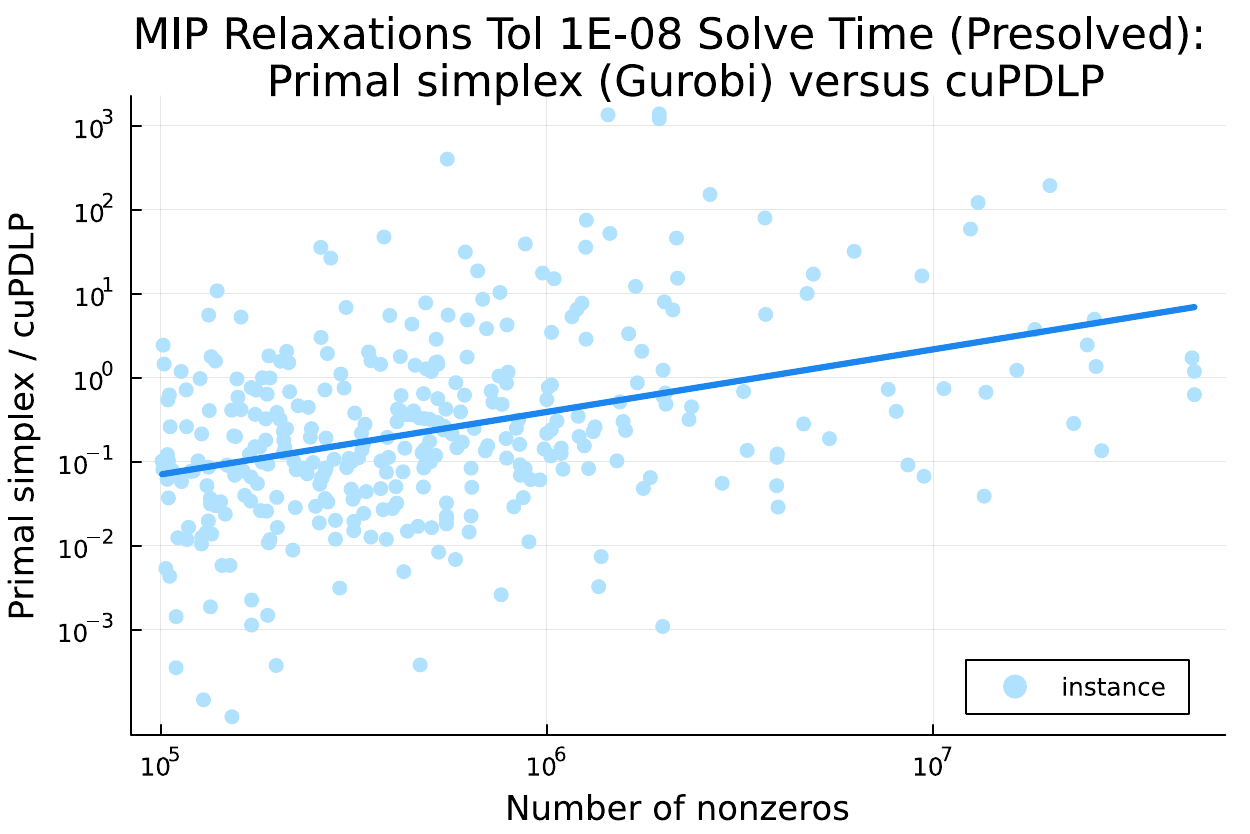}
        & \includegraphics[width=0.33\textwidth]{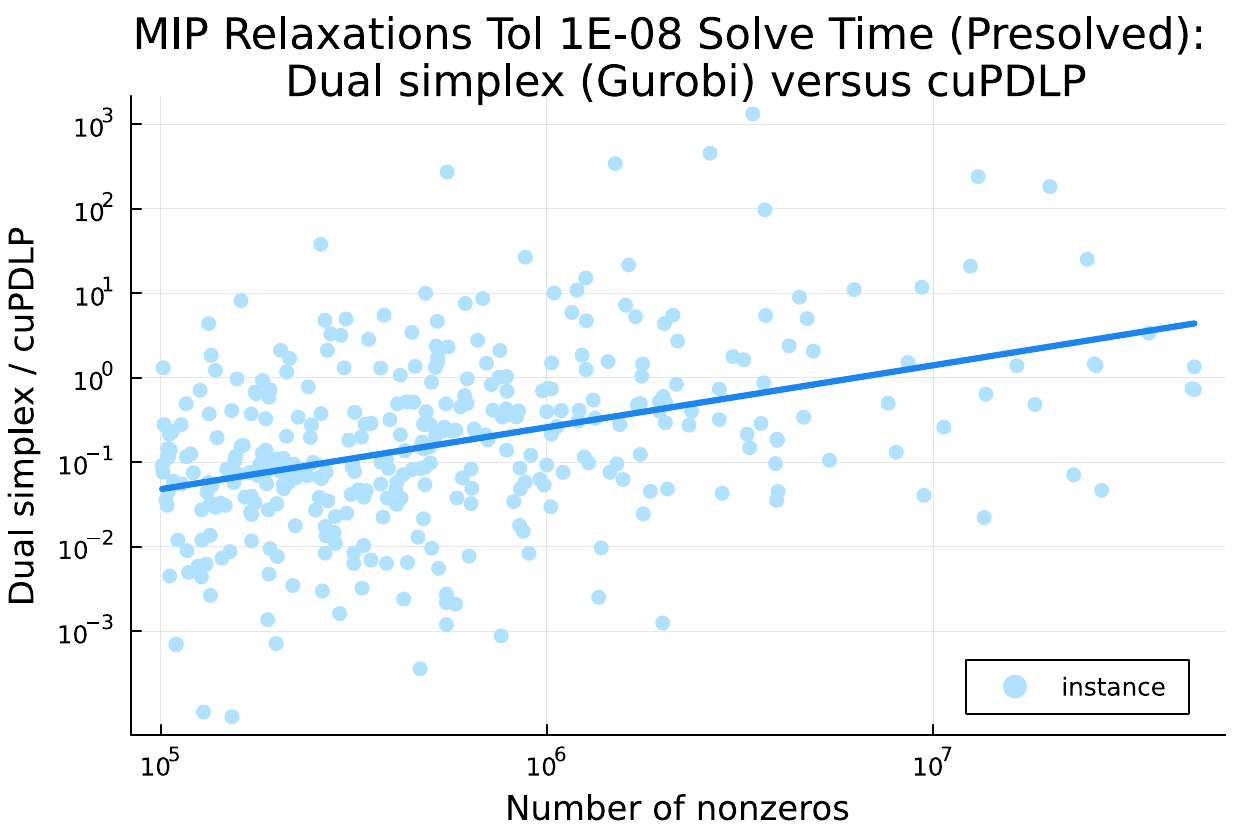}
        & \includegraphics[width=0.33\textwidth]{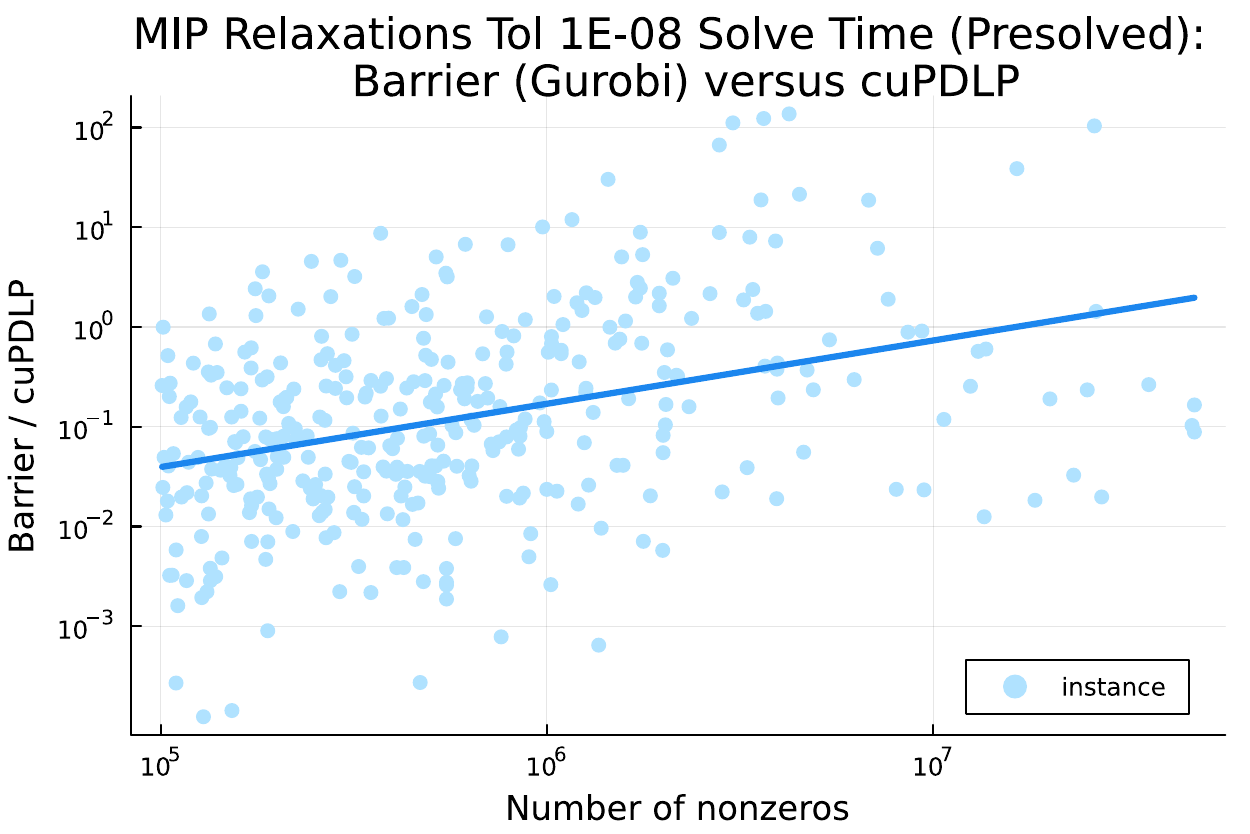}
	\end{tabular}
	\caption{Ratio of Gurobi solve time over cuPDLP.jl solve time for moderate accuracy (top) and high accuracy (bottom). Presolve is used.}
	\label{fig:scatterplot-gurobi-with-presolve}
\end{figure}

Figure \ref{fig:scatterplot-gurobi-no-presolve} and Figure \ref{fig:scatterplot-gurobi-with-presolve} visualize the running time comparison between Gurobi and cuPDLP.jl over the size of the instances with scatter plots. Specifically, each dot in a scatter plot is an instance both methods can solve within the time limit. The x-axes are the number of nonzeros of a \texttt{MIP Relaxations} instance, and the y-axes are the running time ratio of Gurobi (primal simplex, dual simplex, and barrier, respectively) over cuPDLP.jl. Though with higher variance, we can observe a strong positive correlation between the speed-up of cuPDLP.jl against Gurobi versus the number of nonzeros of the problems. We can also observe quite a few instances where cuPDLP is more than 100 times faster than Gurobi different methods.

\begin{figure}[ht!]
	\begin{tabular}{c c c}
		& \includegraphics[width=0.45\textwidth]{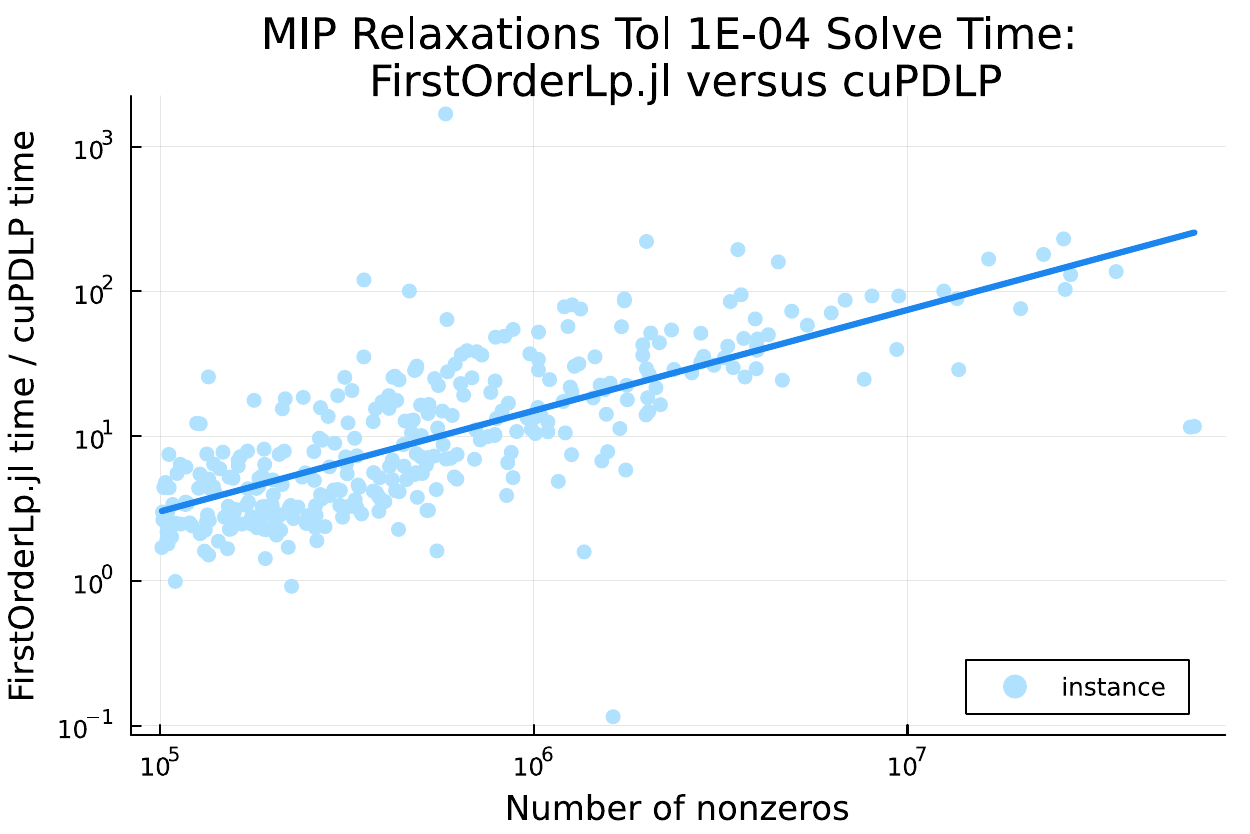}
        & \includegraphics[width=0.45\textwidth]{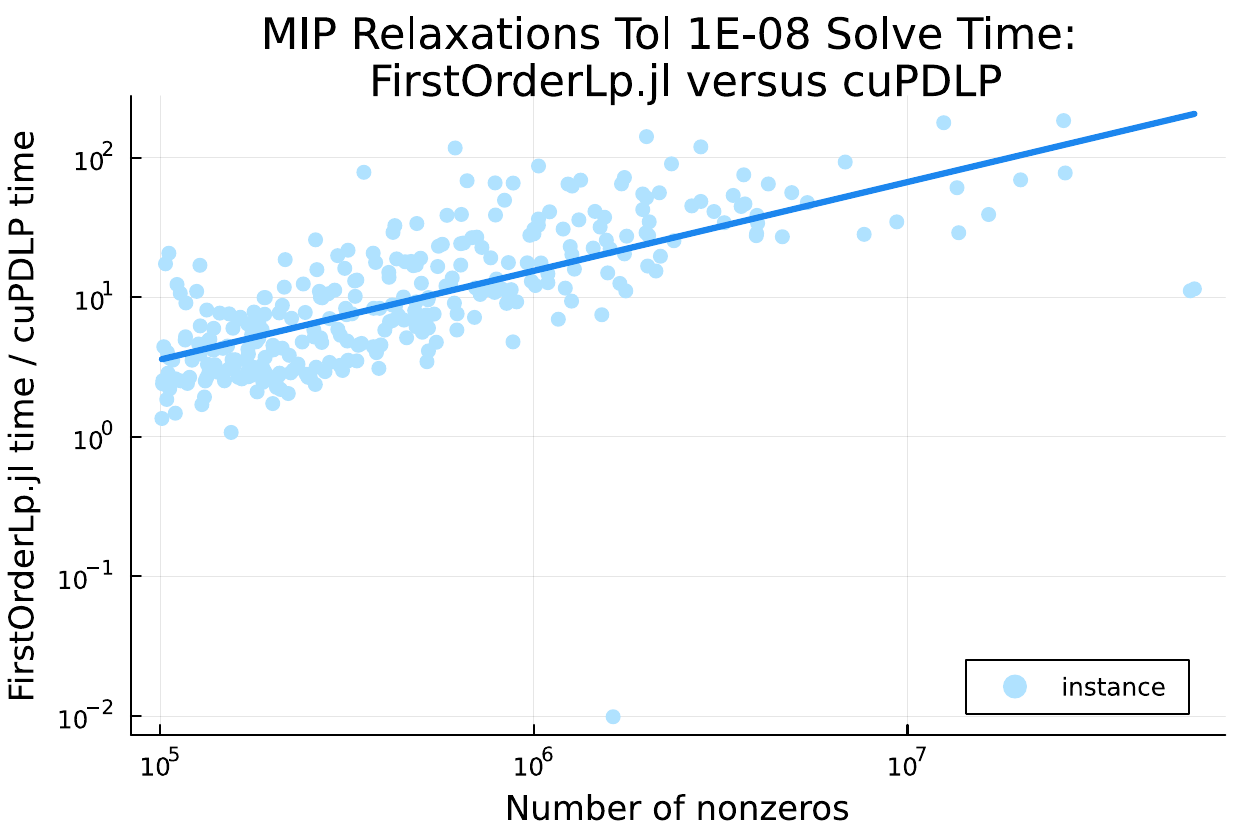}\\
        & \includegraphics[width=0.45\textwidth]{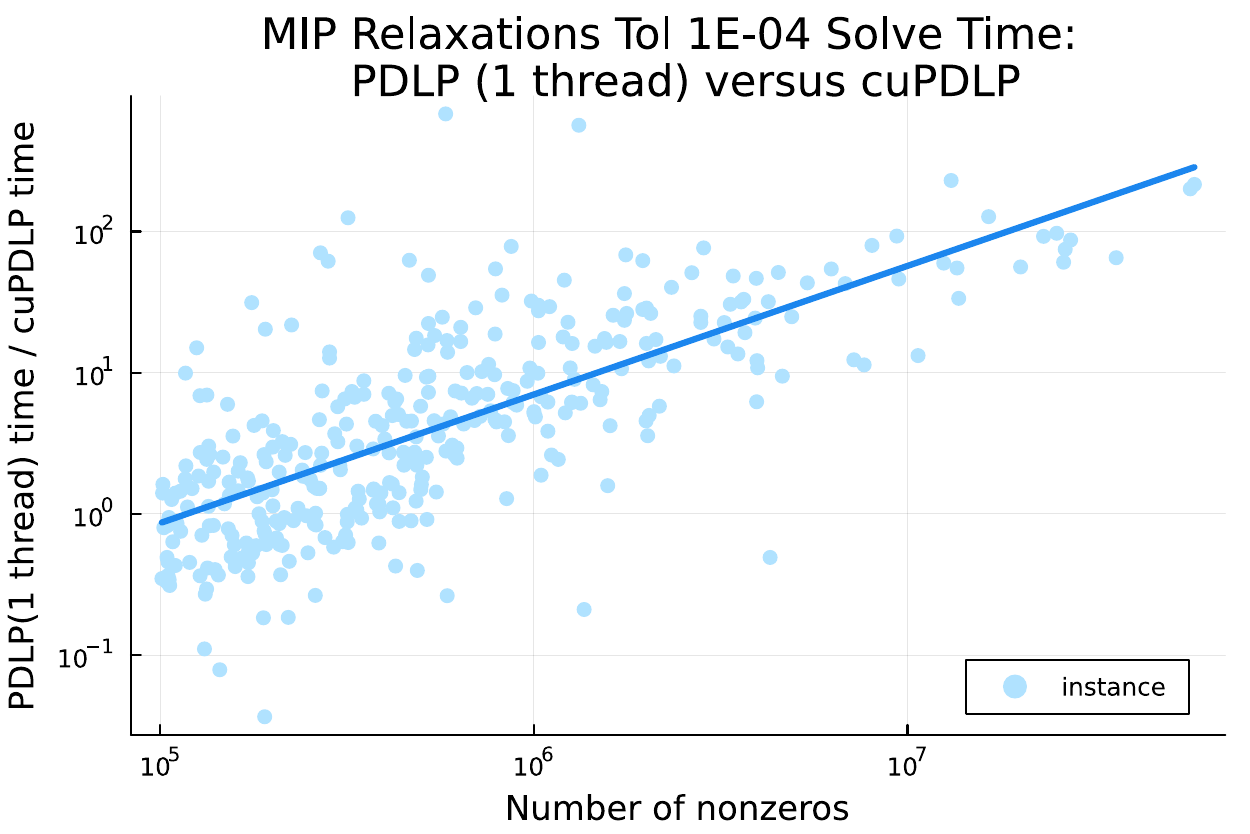}
        & \includegraphics[width=0.45\textwidth]{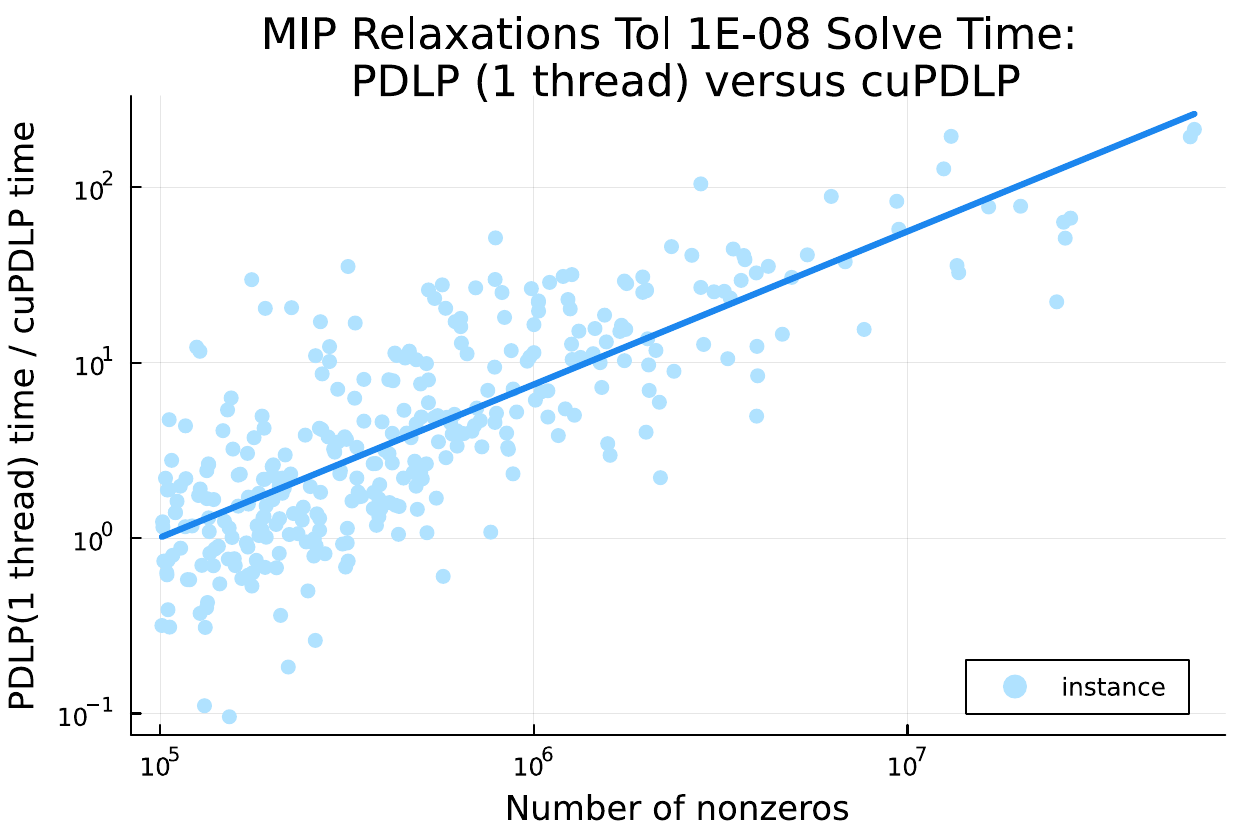}\\
        & \includegraphics[width=0.45\textwidth]{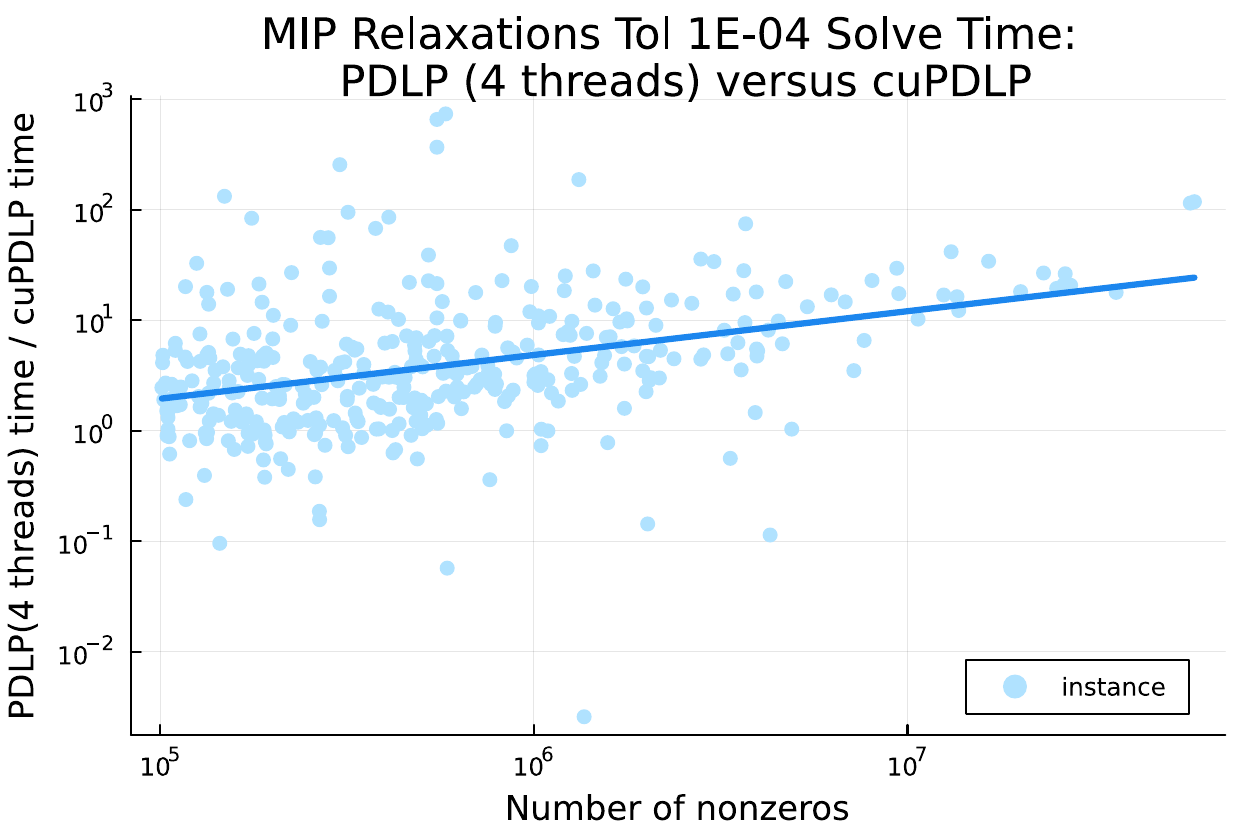}
        & \includegraphics[width=0.45\textwidth]{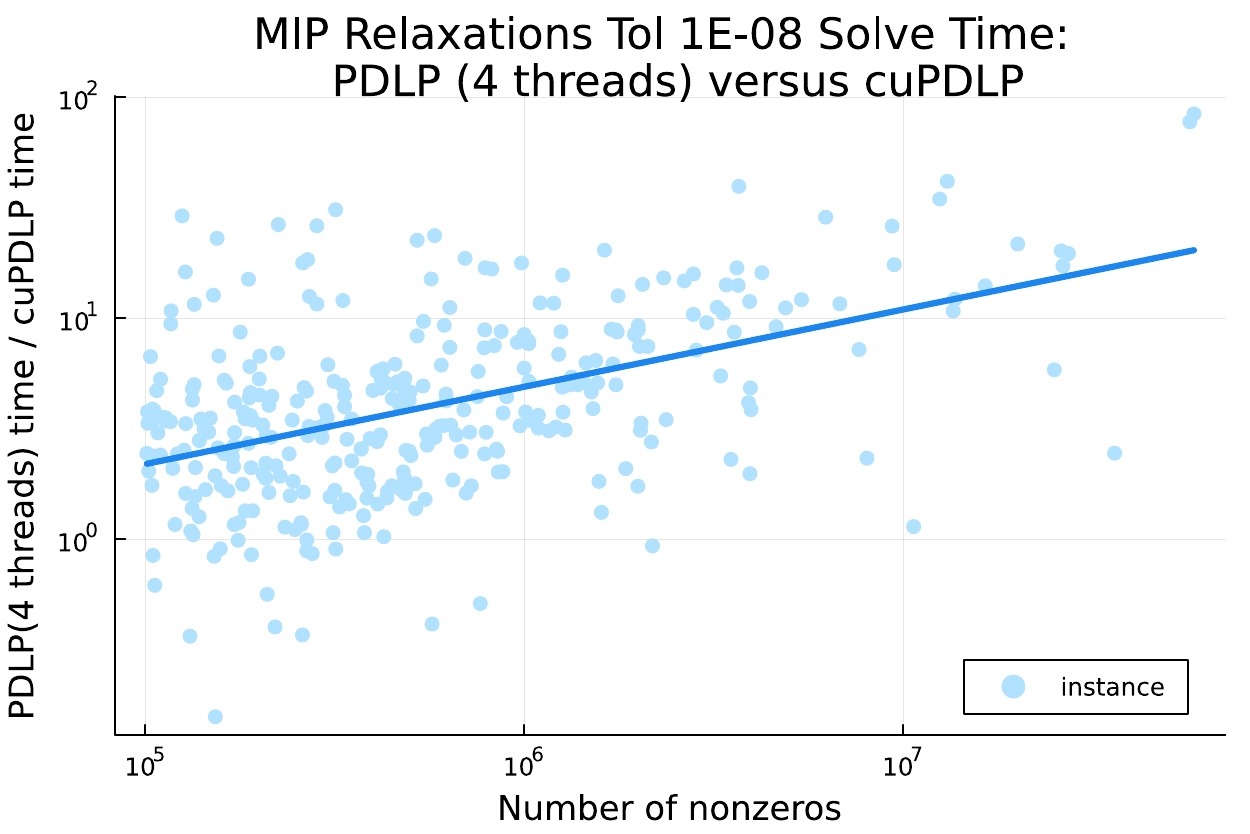}\\
        & \includegraphics[width=0.45\textwidth]{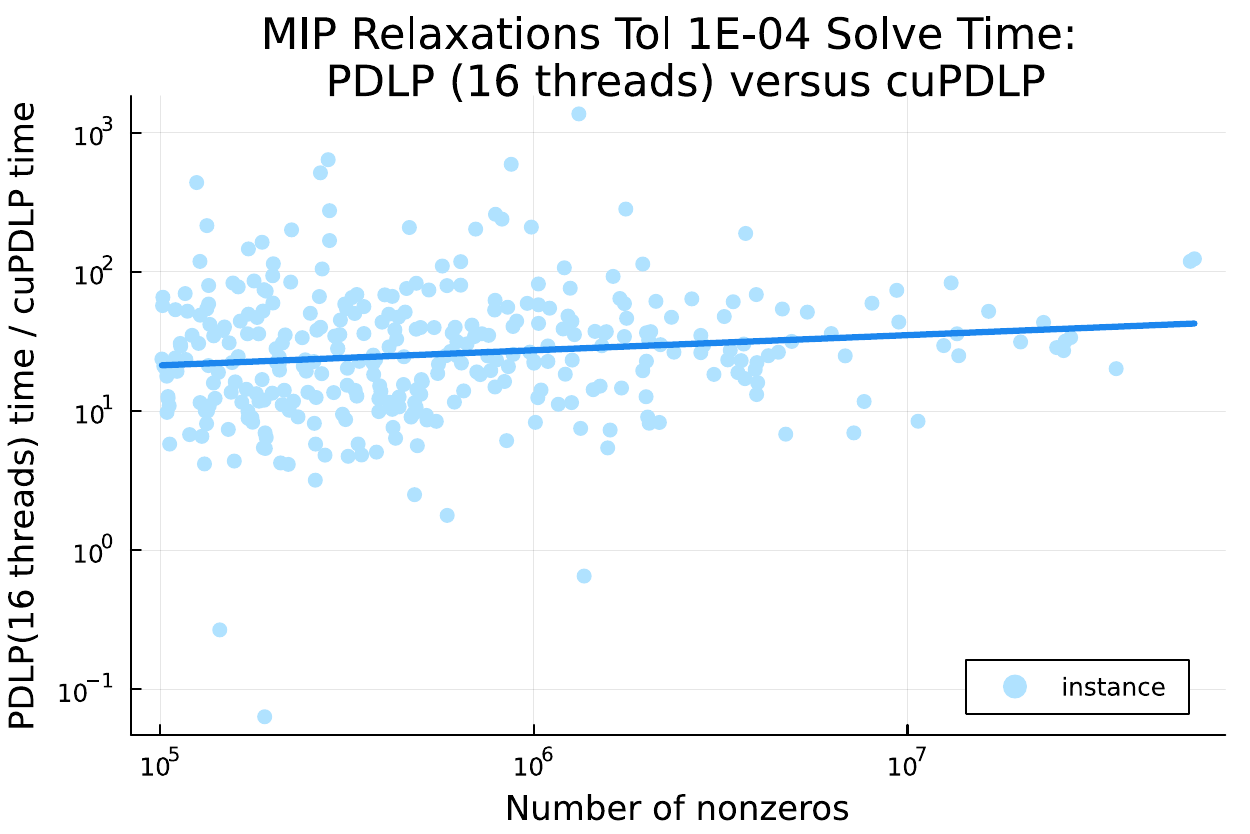}
        & \includegraphics[width=0.45\textwidth]{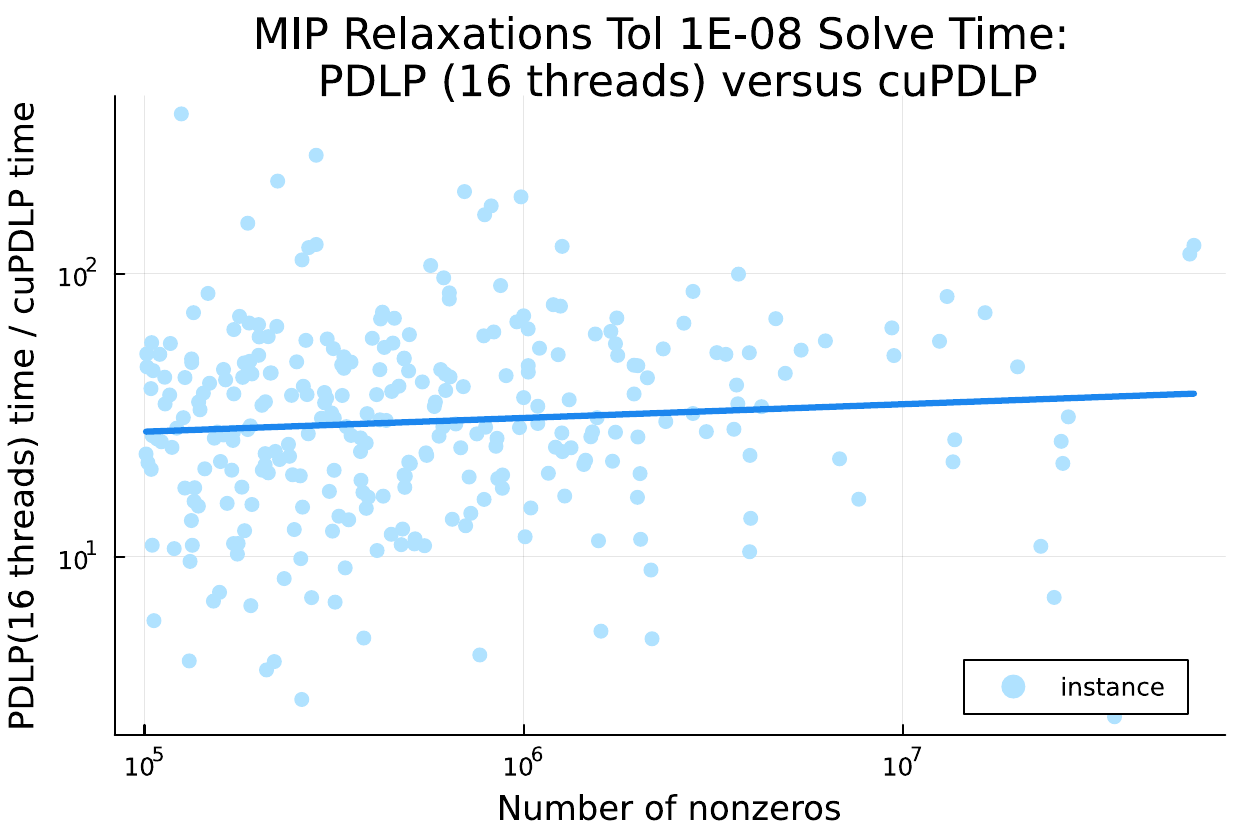}       
	\end{tabular}
	\caption{Ratio of PDLP solve time over cuPDLP.jl solve time for moderate accuracy (top) and high accuracy (bottom).}
	\label{fig:scatterplot-pdlp}
\end{figure}

Furthermore, we visualize the comparison of solve time of PDLP with cuPDLP.jl in Figure \ref{fig:scatterplot-pdlp}. The y-axes are the ratio of FirstOrderLp.jl/PDLP with single thread/PDLP with 4 threads/PDLP with 16 threads solve time and cuPDLP.jl solve time, which represents the speed-up gained by cuPDLP.jl, while the x-axes is the number of nonzeros in each constraint matrix. The increasing trends in the figures imply that GPU-based cuPDLP.jl can gain more significant speed-up over CPU-based counterpart PDLP as we solve larger instances.

\end{document}